\documentclass[onefignum,onetabnum]{siamart171218}

\usepackage{tikz}

\usepackage[title]{appendix}
\usepackage[showframe=false]{geometry}
\usepackage{changepage}

\usepackage{ulem}

\usetikzlibrary{arrows,shapes, positioning, matrix, patterns, decorations.pathmorphing, arrows.meta,automata}

\usepackage{color}
\usepackage{colortbl}
\usepackage{soul}

\usepackage{array}
\hypersetup{colorlinks=true,linkcolor=blue,citecolor=blue,urlcolor=blue}
\usepackage{amsfonts}
\usepackage{amssymb}
\usepackage{mathrsfs}
\usepackage{amsmath}
\usepackage{adjustbox}

\renewcommand{\theequation}{\thesection.\arabic{equation}}

\newcommand{\R}{\mbox{$\mathbb{R}$}}

\newcommand{\C}{\mbox{$\mathbb C$}}

\newcommand{\etal}{{\it et al.}}

\newcommand{\END}{\hfill$\Diamond$}

\newsiamthm{header}{}
\newsiamthm{prop}{Proposition}
\newsiamthm{thm}{Theorem}
\newsiamthm{cor}{Corollary}
\newsiamthm{con}{Conjecture}
\newsiamremark{Def}{Definition}
\newsiamremark{exam}{Example}
\newsiamremark{exams}{Examples}
\newsiamremark{exer}{Exercise}
\newsiamremark{exers}{Exercises}
\newsiamremark{rem}{Remark}
\newsiamremark{rems}{Remarks}

\makeatletter
\@ifundefined{refstepcounter@optarg}{}{%
  \def\refstepcounter@optarg[#1]#2{\refstepcounter{#2}}}
\makeatother

\title{Zero-Hopf networks}
\author{Manuela Aguiar\thanks{Faculdade de Economia, Universidade do Porto, Centro de Matem\'atica da Universidade do Porto, Rua Dr Roberto Frias, 4200-464 Porto, Portugal (\email{maguiar@fep.up.pt}).}
\and Ana Dias\thanks{Dep. Matem\'atica, Faculdade de Ci\^encias, Universidade do Porto, Centro de Matem\'atica da Universidade do Porto, Rua do Campo Alegre, 687, 4169-007 Porto, Portugal (\email{apdias@fc.up.pt}).}
\and Fahimeh Mokhtari\thanks{Vrije Universiteit Amsterdam, Department of Mathematics, NIEUW UNIVERSITEITSGEBOUW, De Boelelaan 1111, 1081 HV Amsterdam, Holland (\email{f.mokhtari@vu.nl}).}}

\headers{Zero-Hopf networks}{M. Aguiar, A. Dias, and F. Mokhtari}

\begin{document}

\maketitle

\begin{abstract}
A zero-Hopf bifurcation can occur for dynamical systems that have a singularity where the linear part has a zero eigenvalue and a pair of purely imaginary eigenvalues.  We  call a zero-Hopf network a regular network whose adjacency matrix has a zero eigenvalue and a pair of complex eigenvalues with zero real part.  In this work, we prove that there is an infinite number of zero-Hopf networks  for which the zero-Hopf bifurcation can occur as a generic codimension-one phenomenon for the class of the network admissible smooth vector fields.   Our results are achieved using the network adjacency matrices and the associated eigenvector structures. General results are presented for four-cell zero-Hopf networks concerning the network set of synchrony spaces. As a case study, we present a classification of all connected zero-Hopf 4-cell regular networks of valency two showing that, up to permutation of the cells,  there are 9 distinct such networks. Moreover, the normal forms of these networks are computed using the method developed by Rink and Sanders (Coupled cell networks: semigroups, Lie algebras and normal forms. {\it Transactions of the American Mathematical Society} {\bf 367} (2015) (5) 3509--3548) for coupled cell networks. It is remarked that these normal forms are less simplified as the input symmetry of the network admissible smooth vector fields restricts the admissible normal form transformations. Nevertheless, the networks share a common coupling function normal form. 
\end{abstract}

\begin{keywords}
Zero-Hopf singularity, steady-state bifurcation, Hopf bifurcation, coupled cell network, zero-Hopf network
\end{keywords}

\begin{AMS}
Primary: 34C23, 37G10, 34C25; Secondary: 34C15
\end{AMS}
 

\section{Introduction}

In the class of smooth vector fields, the zero-Hopf bifurcation can occur  for dynamical systems that have a singularity where the linear part has a zero eigenvalue and a pair of purely imaginary eigenvalues, and generically,  this is  a codimension-two bifurcation. Besides limit cycles and saddle-node bifurcations of periodic orbits, this bifurcation can unfold chaos via Shilnikov homoclinic orbits or even heteroclinic connections (Bykov cycles), see Baldom\'a \etal ~\cite{BIS20} and Dumortier \etal ~\cite{DIKC13}.
In fact, it is the simplest local bifurcation that can create strange attractors and chaotic dynamics (unfolding through Shilnikov homoclinic orbits). 
The  zero-Hopf bifurcation has been studied in application models, such as, the Chua circuit, ~\cite{EL17}, and the FitzHugh–Nagumo system, ~\cite{ELV14}.

One of the contributions of this work is to show that in the set of smooth coupled cell systems, that have structure consistent with a network,  
 the zero-Hopf bifurcation can occur as a generic codimension-one phenomenon.

Our work focusses on the formalisms of Golubitsky, Stewart and collaborators presented in \cite{SGP03, GST05} and Field~\cite{F05}. The structure of the {\it coupled cell systems}  (here, systems of ordinary differential equations) is determined by {\it coupled cell networks} (graphs) with nodes indicating the {\it cells} and edges the {\it couplings} between the cells. We consider  {\it regular networks}, where there is only one kind of cell and coupling, and the number of inputs directed to each cell  is constant. 
This constant is called the  {\it network valency}.

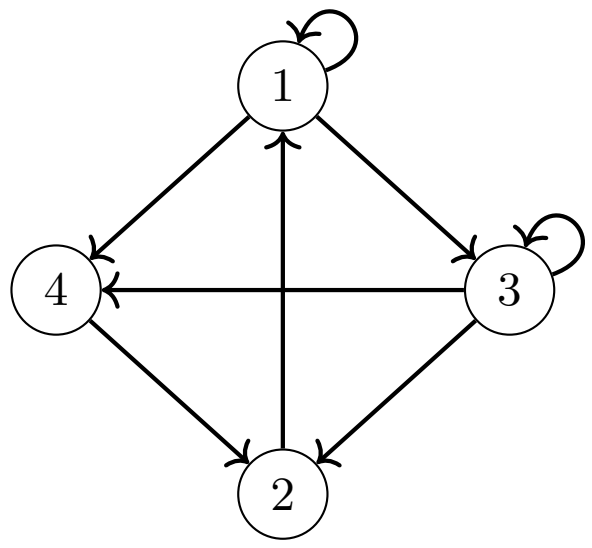
\begin{figure}
\begin{tabular}{c}
\begin{tikzpicture}
 [scale=.15,auto=left, node distance=1.5cm]
 \node[fill=white,style={circle,draw}] (n4) at (4,0) {\small{4}};
  \node[fill=white,style={circle,draw}] (n3) at (24,0) {\small{3}};
 \node[fill=white,style={circle,draw}] (n1) at (14,9)  {\small{1}};
 \node[fill=white,style={circle,draw}] (n2) at (14,-9)  {\small{2}};
 \draw[->, thick] (n1) edge[thick]  node  [above=1pt] {} (n4); 
 \draw[->, thick] (n1) edge[thick]  node  [above=1pt] {} (n3); 
 \draw[->, thick] (n4) edge[thick] node  [above=1pt] {} (n2); 
\draw[->, thick] (n3) edge[thick] node  [below=1pt]  {} (n2); 
\draw[->, thick] (n2) edge[thick] node  [below=1pt]  {} (n1); 
\draw[->, thick] (n3) edge[thick] node  [below=1pt]  {} (n4); 
\draw[->, thick] (n1) to  [in=70,out=20,looseness=5] (n1);
\draw[->, thick] (n3) to  [in=70,out=20,looseness=5] (n3);
\end{tikzpicture}  \\
 \\
 $G$
 \end{tabular}
\caption{Figure 2 of \cite{GL09}. A four-cell regular network $G$ of valency two: all cells and all edges are of the same type, and every cell receives two inputs.  The network adjacency matrix has simple eigenvalues: $2,\, 0,\, \overline{+}i$.}
\label{f:um}
\end{figure}

Golubitsky and Lauterbach~\cite{GL09} have shown that a codimension-one bifurcation from a synchronous equilibrium in a regular network is at linear level isomorphic to a generalized eigenspace of the adjacency matrix of the network, when the dimension of the internal dynamics of each node is greater than one.  As it was remarked in ~\cite{GL09}, in case the cells have  one-dimensional phase space, there are examples where this result does not hold. For those examples, the generic structure of the critical eigenspaces differs from the ones that occur when the cells are assumed to have higher-dimensional phase space, with implications in the analysis of the possible codimension-one steady-state and Hopf synchrony-breaking bifurcations. Our work is precisely motivated by the network example  in~\cite{GL09} pictured in Figure~\ref{f:um}, which we analyse throughout  Section~\ref{sec:Hopfzero}.  

In this work we examine regular network structures that may lead to a codimension-one steady-state and Hopf synchrony-breaking bifurcation, a  {\it zero-Hopf synchrony-breaking bifurcation}. We consider coupled cell systems with one-di\-men\-sional cell phase spaces  associated with the regular networks and whose  critical space, which  is known to be  determined by the Jacobian matrix at the bifurcation point, has the simple zero eigenvalue  and a unique pair of pure imaginary simple eigenvalues. For that, we consider regular networks whose adjacency matrix has a zero eigenvalue and a pair of imaginary eigenvalues, which we  call {\it zero-Hopf networks}.

We take regular networks with four cells, as four is the smallest number of cells to obtain an adjacency matrix 
with four simple eigenvalues: the (nonzero) network valency $k$, scaled so that the network adjacency matrix has a zero eigenvalue and a pair of conjugate imaginary eigenvalues. The type of phenomenon we are addressing here may occur in general regular networks with $n$ cells where $n>4$, in case, the eigenvalues of the network adjacency matrix include the valency $k$, simple eigenvalues $\mu, \mu \pm \omega i$ such that $\mu, \omega \in \R,\,  \mu \not= k,\, \omega \not=0$, and where   the eigenspace sum $E_k \oplus E_{\mu} \oplus E_{\mu \pm \omega i}$ 
corresponds to (or it is contained in) a network synchrony space.

We also investigate normal forms under near-identity coordinate
transformations that preserve the network architecture and input
symmetry \cite{R15}.

For the network considered in our main normal form analysis, the motivating example  from ~\cite{GL09} presented in Figure~\ref{f:um}, we prove in Section~\ref{sec_normal_form} that, at each homogeneous degree $m\geq2$,
the admissible generators whose bracket with the linear part
remains admissible form the one-dimensional space spanned by
\[
\mathbf{M}_m=\sum_{j=1}^{4}x_j^m\frac{\partial}{\partial x_j}.
\]
These generators correspond to applying the same scalar
coordinate change to every cell. Their Lie brackets with
arbitrary admissible vector fields remain admissible, allowing
a normal form procedure that preserves the network structure
at every order. However, the restricted homological image is
one-dimensional at each degree: one can eliminate either the
coefficient of $v^m+w^m$ or that of $u^{m-1}(v+w)$ in the
coupling function, but generally not both. The resulting normal
form therefore retains substantially more nonlinear terms than
a normal form allowing unrestricted coordinate changes.
This identifies the extent of simplification available within
the chosen network-preserving class. Moreover, the action of
these diagonal generators on coupling functions is independent
of the input maps, yielding a common coupling-function
normal form for the additional zero-Hopf 4-cell regular networks of valency two considered here. 

The paper is organized in the following way. In Section~\ref{sec:Hopfzero} we introduce the main definitions that are used throughout the paper, including  regular networks, network admissible coupled cell systems, zero-Hopf networks, network synchrony spaces  and synchrony-breaking zero-Hopf bifurcation. The network in~\cite{GL09} represented in Figure ~\ref{f:um} is used to illustrate all the concepts. 
Section~\ref{sec:zHfourcells} characterizes the possible synchrony spaces of four-cell zero-Hopf networks. Proposition~\ref{prop_possible_sl} shows that there are only four possible distinct synchrony lattice structures. As a corollary, we obtain in Theorem~\ref{thm:infinite} that there is an infinite number of four-cell regular zero-Hopf networks. In Section~\ref{sec_quotients}, the characterization of structural properties of the possible quotient networks of four-cell regular zero-Hopf networks with nontrivial synchrony spaces are obtained. Using that, we obtain Propositions~\ref{prop_steady_state}-\ref{prop_Hopf_bifurcation} concerning steady-state and Hopf bifurcations that can occur generically for those networks. 

The general results obtained in the previous sections are applied in Section~\ref{sec:val2_zHfourcells} for enumerating, up to permutation of the cells, of the zero-Hopf networks of four cells and valency two. Theorems~\ref{thm:um}-\ref{thm:redes_L3} show that, up to permutation of cells, there are only 9 connected such networks. In Section~\ref{sec_normal_form} the normal forms of these networks are computed using the method developed by Rink and Sanders~\cite{R15} for coupled cell networks, Theorem~\ref{thm:first-level-normal-form} and Theorem~\ref{thm:universal-diagonal-normalization}. It is remarked that these normal forms are less simplified in the sense that substantially more nonlinear terms are retained than a normal form allowing unrestricted coordinate change, since the input symmetry of the network admissible smooth vector fields restricts the admissible normal form transformations.

\section{Definitions and motivation example}\label{sec:Hopfzero}

As mentioned in the introduction, our motivational example is the network in~\cite{GL09} represented in Figure ~\ref{f:um}. Throughout this section, we will use this network to illustrate, in addition to the definitions presented, that the type of bifurcation phenomenon we are considering, zero-Hopf bifurcation, can have codimension-one.

We consider {\it regular networks} with constant {\it valency} -- there is only one kind of cell and coupling and the sum of the inputs directed to each cell is constant.  Also, we are only interested in networks that are  {\it connected} -- ignoring the orientation of the connections, there is a path of connections from each cell in the network to any other cell in the network.

\subsubsection{Network adjacency matrix} 
The {\it adjacency matrix} of $G$ denoted by $A_G$ is a $n \times n$ matrix with nonnegative integers, where the $ji$ entry is the number of  edges directed from cell $i$ to cell $j$.
  The adjacency matrix of a regular network of valency $k$ has constant row sum $k$. 

\begin{exam}
The connected network $G$ in Figure~\ref{f:um} has one type of cell and every cell receives two inputs of the same type. Thus, $G$ is a regular four-cell network of valency two. The  network adjacency matrix  is the $4 \times 4$ matrix 
\begin{equation} \label{ex_1_adj_matrix}
A_G = 
\left( 
\begin{array}{cccc}
1 & 1 & 0  & 0 \\
0 & 0 & 1  & 1 \\
1 & 0 & 1 & 0\\
1 & 0 & 1 & 0
\end{array}
\right)\, .
\end{equation} 
\END 
\end{exam}

\subsubsection{Coupled cell systems} 
In the formalisms \cite{SGP03, GST05}, \cite{F05}, 
if $G$ is an $n$-cell regular network of valency $k$,   the coupled cell systems $\dot{x} = F(x,\lambda)$ 
consistent with $G$,  have cells equations  determined by the same  function $f$ with arguments 
determined by the set of cells connected to each cell. More precisely, the evolution of cell $j$ is governed by an equation of the type 
\begin{equation}\label{ccs:integers} 
\dot{x}_j = f(x_j; \overline{x_{j_1}, \ldots, x_{j_k}};\lambda) \quad (j=1,\ldots,n)
\end{equation} 
where $k$ is the network valency, the number of directed edges (of the same type) to cell $j$, which are from cells $j_1, \ldots, j_k$. The bifurcation parameter is $\lambda$.
The overbar indicates the invariance of the function $f$ under permutation of the variables $x_{j_i}$, $i=1,\ldots,k$ for all $\lambda \in \R$. This invariance reflects the fact that all the couplings in the network $G$ have the same type. 
  We assume that cells have cell phase space $\R$ and so $f:\, \R^{k+1} \times \R \to \R$ which is assumed to be smooth. We make the assumption that there is a branch of trivial steady-state solutions where all cells are synchronized for all values of $\lambda$: 
$$
f(0, \lambda) \equiv 0\, .
$$

\begin{exam}Coupled cell systems consistent with the network  $G$ in Figure~\ref{f:um},
have the form:
\begin{equation} \label{e:um}
\begin{array}{l}
\dot{x}_1 = f(x_1; \overline{x_1,x_2}; \lambda)\\
\dot{x}_2 = f(x_2; \overline{x_3,x_4}; \lambda)\\
\dot{x}_3 = f(x_3; \overline{x_1,x_3}; \lambda)\\
\dot{x}_4 = f(x_4; \overline{x_1,x_3}; \lambda)
\end{array}\, .
\end{equation}
Here $f:\, \R^3\times \R \to \R$ is a  one-parameter family of smooth functions, with cells having  phase space $\R$. 
The overbar indicates that $f(u;v,w;\lambda) = f(u;w,v;\lambda)$ for all $u,v,w \in \R$ and $\lambda \in \R$.
\END 
\end{exam}

 If $\alpha$ and $\beta$ are the partial derivatives of $f(u; v_1, \ldots, v_k; \lambda)$  with respect to the first and second variables, respectively, evaluated at the origin, 
 $$
\alpha = \left( \frac{\partial f}{\partial u}\right)_{0}, \quad \beta =  \left(\frac{\partial f}{\partial v_1}\right)_{0}, 
$$
the invariance of $f(u; v_1, \ldots, v_k; \lambda)$ under permutation of the $k$ coordinates $v_1, \ldots, v_k$ implies that 
$$
\beta = \left( \frac{\partial f}{\partial v_i}\right)_0, \quad \forall i \in \{ 1,\ldots, k\}\, .
$$
It follows so that the Jacobian of the admissible coupled cell system for $G$ determined by $f$ as in (\ref{ccs:integers}) 
at the origin  is 
\begin{equation} \label{gen_Jacobian_integer} 
(dF)_0 = \alpha \mathrm{Id}_n + \beta A_G\, .
\end{equation} 
The eigenvalues of $(dF)_0$ are $\alpha + \mu \beta$ where $\mu$ runs through the eigenvalues of $A_G$. Moreover, the eigenvector structure 
for $A_G$ and $(dF)_0$ coincide. 

\begin{exam}
For the network  $G$ in Figure~\ref{f:um}, the eigenvalues of $A_G$ are the valency $2$ and $0, \pm i$. Moreover, the Jacobian matrix of the $G$-admissible systems (\ref{e:um}) at the origin $(x,\lambda) = (0,0)$ has the following form:
$$
(dF)_0 = \alpha \mathrm{Id}_4 + \beta A_G
$$
where $\alpha$ and $\beta$ are the partial derivatives of $f(u;v,w; \lambda)$ with respect to the first and second variables, respectively, evaluated at the origin. The set of eigenvalues of $(dF)_0$ is so  
$$
\left\{ 
\alpha + \mu \beta:\, \mu \in \{ 2, 0, \pm i\} 
\right\} 
= 
\left\{ 
\alpha + 2\beta, \alpha, \alpha \pm i \beta
\right\}, 
$$
and  the Jacobian matrix $(dF)_0$ has  the network adjacency matrix $A_G$  eigenvector structure: 
$$
\hspace{-1mm} 
E_{2}  =  <(1,1,1,1)>, \ 
E_{0}  =   <(1,-1,-1,1)> ,  \  
E_{\pm i}  =  <(-1,0,1,1),\, (1,-2,0,0)> \, .
$$
\END 
\end{exam}

\subsubsection{Zero-Hopf networks}
We have that, if $0$ and $\pm \omega i$, where $\omega \not=0$ are eigenvalues of $A_G$, 
then $\alpha$ and $\alpha \pm \beta \omega i$ are eigenvalues of $(dF)_0$.  

 Assume that:
 \begin{equation}   \label{eq:simple_eig}
 \begin{array}{l} 
  0, \pm \omega i  \mbox{ are simple eigenvalues of $A_G$;}\\
   \mbox{no other eigenvalues of $A_G$ are on the imaginary axis.}
   \end{array} 
  \end{equation} 

 For a codimension-one steady-state and Hopf 
synchrony-breaking bifurcation to occur at the the origin, we assume the 
following conditions on the partial derivatives of $f$ at the origin:  
\begin{equation}   \label{eq:alfa}
\begin{array}{l} 
\alpha = 0 \quad  \mbox{(The bifurcation condition)} \\
 \beta \not= 0  \quad \mbox{(The (generic) nondegeneracy condition)} 
 \end{array} 
 \end{equation}
 As 
$$
\alpha = 0  \Longleftrightarrow \mathrm{Re} (\alpha \pm \beta \omega i )= 0, 
$$
under the conditions   (\ref{eq:simple_eig})-(\ref{eq:alfa}), 
the eigenvalues of $(dF)_0$ include $\alpha + 0 \beta = 0, \alpha \pm \omega i  \beta = \pm \omega i \beta$ and no other eigenvalues of $(dF)_0$ are on the imaginary axis. 
We have so that  the critical subspace of $(dF)_0$  is the  three-dimensional space 
\begin{equation} \label{eq:ogeneralplus}
E_{0}  \oplus E_{\pm \omega  i}, 
 \end{equation} 
 corresponding to the sum of the eigenspaces associated with  the eigenvalues $0$ and $\pm \omega i$ of $A_G$, and we say that the origin is a zero-Hopf equilibrium.  
 Here, we use the eigenspace notation for an eigenvalue $\mu$ of $A_G$: 
$$
E_{\mu} = 
\left\{ 
\begin{array}{l} 
\displaystyle \{ x \in \R^n:\,  \left(A_G - \mu \mbox{Id}_{\R^n}\right) x = 0\} \quad (\mbox{if} \ \mu \in \R);\\
\\
\displaystyle  \{ x \in \R^n:\,  \left(A_G - \mu \mbox{Id}_{\R^n}\right)   \left(A_G - \overline{\mu} \, \mbox{Id}_{\R^n}\right) x = 0\} \quad  (\mbox{if} \ \mu \in \C \setminus \R)\, .
\end{array}  
\right. 
$$

In this work we address connected regular networks giving rise to bifurcation problems where the critical space is precisely of the form $E_0 \oplus E_{\pm \omega i}$, where $E_0,\, E_{\pm \omega i}$ are the eigenspaces of $A_G$. We call those networks zero-Hopf networks:

\begin{Def}
A {\it zero-Hopf network} is a regular $n$-cell network $G$ of valency $k>0$ whose adjacency matrix $A_G$ has simple eigenvalues $k,0,\pm \omega i$, where 
$\omega \in \R \setminus \{0\}$,   and no other eigenvalues on the imaginary axis. 
\hfill $\Diamond$
\end{Def}

\begin{exam}
The network in Figure~\ref{f:um} is a zero-Hopf network. 
For a codimension-one steady-state synchrony-breaking bifurcation to occur at the the origin, we assume 
the conditions on the partial derivatives of $f$ at the origin given by  (\ref{eq:alfa}) and 
 the critical subspace of $(dF)_0$  is the  three-dimensional space 
\begin{equation} \label{eq:oplus}
E_{0}  \oplus E_{\pm i}, 
 \end{equation} 
so the origin is a zero-Hopf equilibrium. 
\END
\end{exam}

\subsubsection{Synchrony spaces}

A {\it polydiagonal space} of $\R^n$ is a subspace  defined by certain equalities of coordinates. Such a space is called a 
{\it network synchrony space} when is  left invariant under the flow of every  network admissible coupled cell system.  
Since $G$ is regular, we have that 
the eigenspace associated with the network valency $k$, 
$$
\Delta_0  = <(1,1, \ldots, 1)> = E_k
$$ 
is always one such synchrony space, which is called the {\it full synchrony space}.   

\begin{exam}
For the network in Figure~\ref{f:um}, note that the full diagonal space 
$$\Delta_0 = \{x \in \R^4:\, x_1 = x_2 = x_3 = x_4\}$$ 
is left invariant under the equations (\ref{e:um}), independently of the chosen function $f$. Thus any solution of (\ref{e:um}) with initial condition in $\Delta_0$  
remains there for all time and 
has all cells evolving in time in the same way. 
\END
\end{exam}

A network can have additional synchrony spaces $\Delta_i$. Moreover, 
 the corresponding restriction to each $\Delta_i$ of the network coupled cell systems are again coupled cell systems associated with a network with fewer cells, 
 called the {\it quotient network} of $G$ by $\Delta_i$, that we are denoting by $Q_i$, see \cite{GST05}. It is also said  that the network $G$ is a {\it lift network}
  of that (quotient) network, see \cite{ADGL09}. In general, fixing a $k$-cell network $Q$, there are many $n$-cell networks $G$, where $n>k$, 
  admitting a synchrony space such the quotient of $G$ by that space is $Q$. However, fixing a $n$-cell network $G$ and a synchrony space of $G$, the quotient network obtained from $G$ using that synchrony space is unique. 
 
 The set of synchrony spaces of a network $G$, using the partial order of inclusion, has the structure of a finite lattice, called the {\it network synchrony lattice} and denoted here 
 by  $\mathcal{L}_G$, see Stewart~\cite{S07}. For this lattice, the {\it join} of two synchrony spaces is the smallest synchrony subspace containing them. The {\it meet} of two synchrony spaces  is the largest synchrony subspace contained in both  (the intersection of the two synchrony subspaces). For regular networks, the {\it bottom} of the lattice is always the {\it full synchrony space}  
 $\Delta_0$ and the {\it top} of the lattice is the {\it network phase space} ($\R^n$, in our case).  
 A synchrony space which is not the bottom neither the top of the lattice is said to be {\it nontrivial}.

 One way of describing the synchrony spaces of a network is by finding the polydiagonal spaces that are left invariant under the network adjacency matrix. In fact, it is possible to describe the synchrony lattice of a network using the eigenvector structure of the network adjacency matrix. See Aguiar and Dias~\cite{AD14}. 
See also  Kamei~\cite{K09} in case the adjacency matrix has simple eigenvalues.

\begin{figure}
{\tiny 
\begin{tabular}{ccccc}
\begin{tikzpicture}
 [scale=.15,auto=left, node distance=1.5cm]
 \node[fill=white,style={circle,draw}] (n1) at (4,0) {$\small{\overline{1}}$};
\draw[->, thick] (n1) to  [in=80,out=20,looseness=5] (n1);
\draw[->, thick] (n1) to  [in=60,out=20,looseness=5] (n1);
\end{tikzpicture}  & 
\begin{tikzpicture}
 [scale=.15,auto=left, node distance=1.5cm]
 \node[fill=white,style={circle,draw}] (n1) at (4,0) {$\overline{\small{1}}$};
  \node[fill=white,style={circle,draw}] (n2) at (14,0) {$\overline{\small{2}}$};
 \draw[->, thick] (n1) edge[thick]  node  [above=1pt] {} (n2); 
 \draw[->, thick] (n2) edge[thick]  node  [above=1pt] {} (n1); 
\draw[->, thick] (n1) to  [in=70,out=20,looseness=5] (n1);
\draw[->, thick] (n2) to  [in=20,out=70,looseness=5] (n2);
\end{tikzpicture}  & 
\begin{tikzpicture}
 [scale=.15,auto=left, node distance=1.5cm]
 \node[fill=white,style={circle,draw}] (n3) at (4,0) {$\overline{\small{3}}$};
  \node[fill=white,style={circle,draw}] (n1) at (14,0) {$\overline{\small{1}}$};
   \node[fill=white,style={circle,draw}] (n2) at (9,-9) {$\overline{\small{2}}$};
 \draw[->, thick] (n3) edge[thick]  node  [above=1pt, near end] {{\tiny $2$}} (n2); 
\draw[->, thick] (n3) to  [in=70,out=20,looseness=5] (n3);
\draw[->, thick] (n1) to  [in=20,out=70,looseness=5] (n1);
\draw[->, thick] (n1) edge[thick]  node  [above=1pt] {} (n3); 
\draw[->, thick] (n2) edge[thick]  node  [above=1pt] {} (n1); 
\end{tikzpicture}     &
\begin{tikzpicture}
 [scale=.15,auto=left, node distance=1.5cm]
\node (n3) [ellipse,draw]  at (19,9) {$\Delta_3; 2,\, 0,\, \pm i$};
\node (n2) [ellipse,draw]  at (24,2) {$\Delta_2; 2,\, \pm i$};
\node (n1) [ellipse,draw]  at (10,-3) {$\Delta_1;\, 2,\, 0$}; 
\node (n0) [ellipse,draw] at (19,-9)  {$\Delta_0;\, 2$};
\draw (n0) to (n1);
\draw (n0) to (n2);
\draw (n1) to  (n3);
\draw (n2) to  (n3);
\end{tikzpicture} &\hspace{-10mm} 
$\begin{array}{l}
\Delta_3 = \R^4\\
\Delta_2 = \{ x:\, x_3 = x_4\}\\
\Delta_1 = \{ x:\, x_1 = x_4,\, x_2 = x_3\}\\
\Delta_0 = \{ x:\, x_1 = x_2 = x_3 = x_4\}
\end{array}$\\
$Q_0$ & $Q_1$ & $Q_2$ & $\mathcal{L}_G$ & 
\end{tabular}}
\caption{Three quotient networks $Q_i = G/ \Delta_i$, for $i=0,1,2$ and the synchrony lattice $\mathcal{L}_G$ of the network $G$  in Figure~\ref{f:um}. In the lattice, the synchrony spaces $\Delta_i$ are annotated with the adjacency matrix  eigenvalues of the corresponding quotient network $Q_i$.}
\label{f:dois}
\end{figure}
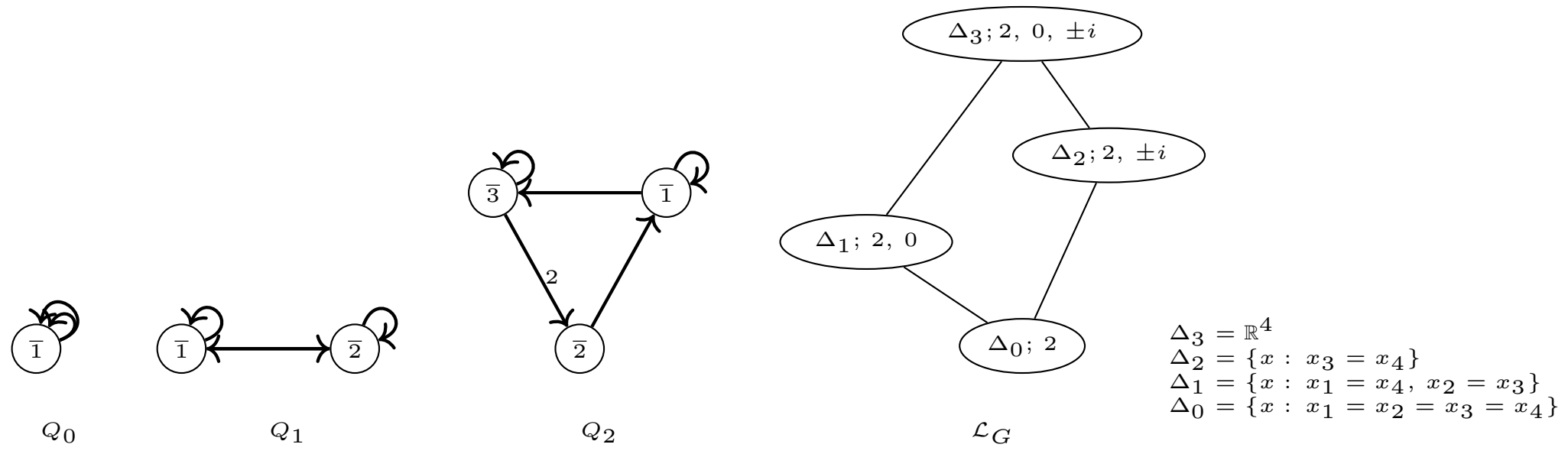

\begin{exam}
The  network in Figure~\ref{f:um} has more synchrony spaces, besides the trivial ones, $\Delta_0$ and $\R^4$. Indeed, 
the  polydiagonal subspaces 
$$
\begin{array}{ll}
\Delta_1  =  \{ x:\, x_1 = x_4,\, x_2 = x_3\}, \quad & \quad \Delta_2  =  \{ x:\, x_3 = x_4\} 
\end{array}
$$
are also synchrony subspaces of $G$ as they are left invariant under the flow of equations (\ref{e:um}). It has so two nontrivial synchrony spaces, the two-dimensional $\Delta_1$ and the three-dimensional $\Delta_2$. Each restriction  of equations (\ref{e:um}) to $\Delta_i$ is a coupled cell system associated with a network with fewer cells, the quotient network of $G$ by $\Delta_i$,   denoted  by $Q_i$. See in Figure~\ref{f:dois} the three quotient networks of $G$ by $\Delta_i$, for $i=0,1,2$. We also say that $G$ is a {\it lift network} of each of the three quotients. We see on the right of Figure~\ref{f:dois} the synchrony lattice of $G$, where the synchrony spaces are annotated with the eigenvalues of the adjacency matrix of the corresponding quotient network. 

 \END
\end{exam}

\subsubsection{Synchrony-breaking zero-Hopf bifurcation}

 Returning to the network example in Figure~\ref{f:um}, we will see that the degeneracy bifurcation condition on $(dF)_0$, recall (\ref{eq:alfa}),  leads to a (codimension-one) zero-Hopf bifurcation. 

 We assume the bifurcation point is $(0,0,0,0)$ and that the bifurcation occurs for $\lambda =0$.  We make the assumption that there is a branch of trivial steady-state solutions where all cells are synchronized for all values of $\lambda$, that is, $f(0, \lambda) \equiv 0$. We ask about the possible bifurcations that will occur if the Jacobian of the system (\ref{e:um}) at $(x,\lambda) = (0,0)$ is singular.  Here, we  assume  the critical space is not contained in the full synchronized space $\Delta_0$, in order that a {\it synchrony-breaking bifurcation} may occur. 

We make the following observation concerning 
the relations of the eigenspaces of $A_G$ ($(dF)_0$) and the nontrivial network synchrony spaces, and recall that $\Delta_0 = E_{2}$: 
$$
\begin{array}{lll}
 \Delta_1 = E_{0} \oplus \Delta_0, \ & \ 
\Delta_2 = E_{\pm i} \oplus \Delta_0 \, .
\end{array}
$$
Thus,  the bifurcation problem for equations (\ref{e:um}) under the assumptions (\ref{eq:alfa}) restricted to $\Delta_1$ has the one-dimensional critical space $E_0$  and corresponds to a standard codimension-one (one-dimensional) steady-state bifurcation problem. Also, the  bifurcation problem for equations (\ref{e:um}) under the assumptions (\ref{eq:alfa}) restricted to $\Delta_2$ has the two-dimensional critical space $E_{\pm i}$, corresponding to a standard codimension-one (two-dimensional) Hopf bifurcation problem. 
Moreover, equations (\ref{e:um}) restricted to $\Delta_1$ 
are consistent with the 2-cell regular (quotient) network $Q_1$ in Figure~\ref{f:dois}  and equations (\ref{e:um}) restricted to $\Delta_2$ are consistent with the 3-cell regular (quotient) network  $Q_2$ in Figure~\ref{f:dois}. We see here in this example that under one degeneracy bifurcation condition on $(dF)_0$, recall (\ref{eq:alfa}), two bifurcations occur. Moreover, small perturbations of the system (\ref{e:um}) preserving the $G$ network structure will give the same bifurcations.

As a simulation example, consider, now, the particular realization of 
 system (\ref{e:um}) where 
\begin{equation} \label{e:rea_um}
f(u; \overline{v, w}; \lambda) = \lambda u -0.000001 u^2 + v + w 
\end{equation}
and note that $\alpha = f_u(0) = 0$ and $\beta= f_v(0) = 1 \not=0$. For the cell dynamics (\ref{e:rea_um}), system (\ref{e:um}) restricted to the synchrony space 
$\Delta_2  =  \{ x \in \R^4:\, x_3 = x_4\}$ is 
\begin{equation} \label{e:simul1}
\begin{array}{l}
\dot{x}_1 =      \lambda x_1 -0.000001 x_1^2 + x_1 + x_2 \\
\dot{x}_2 =      \lambda x_2 -0.000001 x_2^2 + 2x_3 \\
\dot{x}_3 =      \lambda x_3 -0.000001 x_3^2 + x_1 + x_3 
\end{array}\, .
\end{equation}

\begin{figure}[htb]
\vspace{-45mm}
\centerline{
\includegraphics[width=0.7\textwidth]{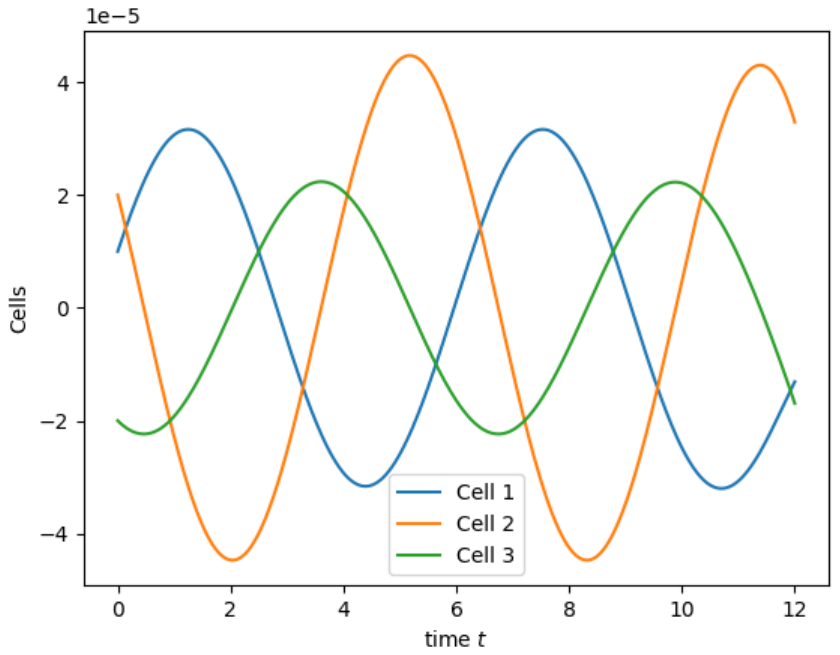}
} 
\vspace{-45mm}
\caption{Time series for a periodic solution of the system (\ref{e:simul1}) obtained by Hopf bifurcation from the origin at $\lambda=0$. 
The lines correspond to the cells $x_1, x_2$ and $x_3$, respectively. The amplitude and phase relations between cells are approximately to 
$|A_2|/|A_1|  =\sqrt{2}, \ |A_3|/|A_1| = \sqrt{2}/2$ and $\phi_2 - \phi_1 = 5/8$, $\phi_3 - \phi_1 = 3/8$, respectively.}
\label{F:0_Hopf_perturbation1}
\end{figure}

The linearization  of (\ref{e:simul1}) at $(0,\lambda)$ has the form 
$$
\left( 
\begin{array}{ccc}
\lambda + 1 & 1 & 0 \\
0 & \lambda  &   2 \\
1 & 0 & 1 + \lambda  
\end{array}
\right) = 
\lambda  \mathrm{Id}_3 + 
\left( 
\begin{array}{ccc}
1 & 1 & 0 \\
0 & 0 &   2 \\
1 & 0 & 1 
\end{array}
\right) 
$$
with eigenvalues $\lambda + 2,  \lambda  \pm i$. At $\lambda = 0$, the eigenvalues are $2, \pm i$. Moreover,  the pair of complex eigenvalues $ \lambda  \pm i$ cross 
the imaginary axis  with nonzero speed at $\lambda=0$. Note that $f_{u \lambda} (0) =1 \not=0$. A standard Hopf bifurcation occurs for the system (\ref{e:simul1}) from the origin at $\lambda =0$, leading to a unique branch of periodic solutions. Moreover, it is known that to the lowest order in $\lambda$, the amplitude of these solutions grow at a rate of  $\lambda^{1/2}$. Observe that the critical eigenvector $(z_1, z_2, z_3) = \left( 1,-1+i, (-1/2)(1+i)\right)$ of the linearization associated with the imaginary eigenvalue $i$, does not depend on the form of the vector field. It is an eigenvector of the matrix that describes the adjacency matrix of the new network, the 3-cell regular (quotient) network  $Q_2$ in Figure~\ref{f:dois}, associated with the coupled cell system restricted to $\Delta_2$. We note that $Q_2$ is the network 24 of Leite and Golubitsky~\cite{LG05}. Following the discussion of \cite[pp 2353-2355]{LG05}, to lowest order in $\lambda$, the ratios between the maximum amplitudes $|A_j|$ of the projections of the periodic solution in cell $j$ are approximately given by $|A_2|/|A_1|  =  |z_2| = \sqrt{2}$ and $|A_3|/|A_1| = |z_3| = \sqrt{2}/2$ and are independent of $\lambda$. Moreover, taking for $j=1,2,3$, $z_j = e^{2 \pi i \phi_j}y_j$ for $y_j >0$ and $0 \leq  \phi_j < 1$, as $z_1=1$ and so $\phi_1 =0$, we obtain the relative phase shifts $\phi_2 - \phi_1$ and $\phi_3 - \phi_1$ between cell coordinate projections of the bifurcating branch of the periodic solutions given by  $\phi_2 -\phi_1 = 3/8,\, \phi_3 - \phi_1 = 5/8$ or $\phi_2 -\phi_1 = 1-3/8 = 5/8,\, \phi_3 - \phi_1 = 1-5/8 = 3/8$. See Figure~\ref{F:0_Hopf_perturbation1}.  Observe this branch of periodic solutions of (\ref{e:simul1}) is  the branch of periodic solutions of the original 4-cell 
coupled cell system (\ref{e:um}) where cell 4 is syncronized with cell 3, that is $x_3(t) = x_4(t)$ for all time. \END 

\section{Synchrony lattice structures of four-cell zero-Hopf networks}\label{sec:zHfourcells}

Assuming one-dimensional cell dynamics and considering synchrony-breaking bifurcations, a zero-Hopf network has at least four cells. Therefore, from now on, we will focus our attention on four-cell zero-Hopf networks. 
In this section we aim to characterize the synchrony subspaces of  four-cell zero-Hopf networks and their lattices.
We start by showing that  four-cell zero-Hopf networks are, at least, valency two. 

\begin{thm}
There is no four-cell zero-Hopf regular network with valency one.
\end{thm}

\begin{proof}
Consider the  adjacency matrix $A=[a_{ij}]_{4\times4}$ of a four-cell regular network with valency one. If the network is zero-Hopf then, besides a pair of pure imaginary eigenvalues, the matrix $A$ has the eigenvalues $1$ and $0$ and, thus, trace equal to $1$. We have, then that the principal diagonal  of $A$ has one entry equal to $1$ and the others are zero. Without loss of generality, we can assume that $a_{11}=1$ and $a_{ii}=0$, $i=2,3,4$. Thus, $a_{1j}=0$, $j=2,3,4$. As $1$ and $0$ are both  simple eigenvalues of $A$, then there must be another entry in the first column of $A$ equals to $1$ and the others $0$. Without loss of generality, we can assume $a_{21}=1$ and $a_{i1}=0$, $i=3,4$. Since $a_{34}=1 - a_{32}$ and $a_{43}=1 - a_{42}$, and since $a_{32}$ and $a_{42}$ can take one of the values $0$ or $1$, we are left with four possible choices for the matrix $A$ where none of them has a pair of pure imaginary eigenvalues. This shows that $A$ cannot be a matrix of a zero-Hopf network.
\end{proof}

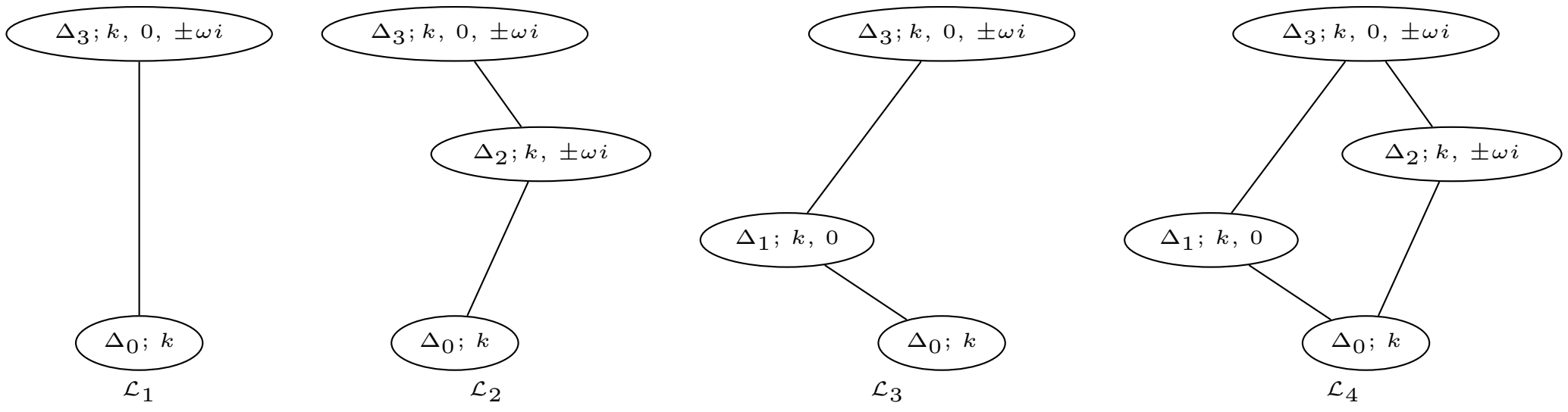
\begin{figure}[!ht]
{\tiny \begin{tabular}{cccc}
\begin{tikzpicture}
 [scale=.15,auto=left, node distance=1.5cm]
\node (n3) [ellipse,draw]  at (19,9) {$\Delta_3; k,\, 0,\, \pm \omega i$};
\node (n0) [ellipse,draw] at (19,-9)  {$\Delta_0;\, k$};
\draw (n0) to (n3);
\end{tikzpicture} &
\begin{tikzpicture}
 [scale=.15,auto=left, node distance=1.5cm]
\node (n3) [ellipse,draw]  at (19,9) {$\Delta_3; k,\, 0,\, \pm \omega  i$};
\node (n2) [ellipse,draw]  at (24,2) {$\Delta_2; k,\, \pm \omega  i$};
\node (n0) [ellipse,draw] at (19,-9)  {$\Delta_0;\, k$};
\draw (n0) to (n2);
\draw (n2) to (n3);
\end{tikzpicture}  &
\begin{tikzpicture}
 [scale=.15,auto=left, node distance=1.5cm]
\node (n3) [ellipse,draw]  at (19,9) {$\Delta_3; k,\, 0,\, \pm \omega i$};
\node (n1) [ellipse,draw]  at (10,-3) {$\Delta_1;\, k,\, 0$}; 
\node (n0) [ellipse,draw] at (19,-9)  {$\Delta_0;\, k$};
\draw (n0) to (n1);
\draw (n1) to (n3);
\end{tikzpicture} & 
\begin{tikzpicture}
 [scale=.15,auto=left, node distance=1.5cm]
\node (n3) [ellipse,draw]  at (19,9) {$\Delta_3; k,\, 0,\, \pm \omega i$};
\node (n2) [ellipse,draw]  at (24,2) {$\Delta_2; k,\, \pm \omega i$};
\node (n1) [ellipse,draw]  at (10,-3) {$\Delta_1;\, k,\, 0$}; 
\node (n0) [ellipse,draw] at (19,-9)  {$\Delta_0;\, k$};
\draw (n0) to (n1);
\draw (n0) to (n2);
\draw (n1) to  (n3);
\draw (n2) to  (n3);
\end{tikzpicture} \\
${\mathcal L}_1$ & ${\mathcal L}_2$ & ${\mathcal L}_3$ & ${\mathcal L}_4$
\end{tabular}}
\caption{The possible synchrony lattice structures for zero-Hopf four-cell regular networks of valency $k$. The synchrony spaces $\Delta_i$ are annotated with the adjacency matrix eigenvalues of the corresponding quotient network: $\Delta_0 = E_k$, $\Delta_1 = E_k \oplus E_0$, $\Delta_2 = E_k \oplus E_{\pm \omega i}$ and $\Delta_3 = \R^4$. 
Here, we assume $k>0$ and $\omega\not=0$.}
\label{f:lista}
\end{figure}

\begin{prop}  \label{prop_possible_sl}
If $G$ is a four-cell zero-Hopf  regular network then 
the synchrony lattice  $\mathcal{L}_G$ of $G$ has  one of the four lattice structures in Figure~\ref{f:lista}. In particular, 
it follows that $G$ can have:\\
(i) no nontrivial synchrony subspaces;\\
(ii) one nontrivial synchrony subspace, two-dimensional or three-dimensional;\\
(iii) two nontrivial synchrony subspace, one two-dimensional and the other three-dimensional.
\end{prop}

\begin{proof} By ~\cite{AD14}, the synchrony subspaces of a regular network $G$ are generated by the eigenvectors, or generalized eigenvectors, of the adjacency matrix $A_G$ of $G$.  More precisely, a synchrony space is a direct sum of invariant spaces of the network adjacency matrix which has the form of a polydiagonal space. 

Let $G$ be a zero-Hopf four-cell regular network with valency $k>1$,  where the adjacency matrix $A_G$ has simple eigenvalues 
$k, 0, \pm \omega i$, with $\omega \not= 0$. It follows so that, the possible  invariant spaces of $A_G$ and corresponding dimensions are:
\begin{center}
\begin{tabular}{|l|c|}
\hline 
$A_G$-invariant space & Dim. \\
\hline 
$\Delta_0 = E_k = \{ (x,x,x,x):\, x \in \R\}$ & $1$ \\
$E_k \oplus E_0$ & $2$ \\
$E_k \oplus E_{\pm \omega i}$ & $3$ \\
$E_0 \oplus E_{\pm \omega i}$ & $3$ \\
$E_k \oplus E_0 \oplus E_{\pm \omega i} =  \R^4$ & $4$\\
\hline 
\end{tabular}
\end{center}
The trivial synchrony subspaces are the full synchony subspace $\Delta_0 = E_k$ and the phase space $E_k \oplus E_0 \oplus E_{\pm \omega i} = \R^4$. Besides these, 
the possible nontrivial synchrony subspaces are given by $E_k \oplus E_{\mu}$, with $\mu = 0$, or $\mu = \pm  i$, in case they are polydiagonals. Note that any 
polydiagonal contains the trivial synchrony space $\Delta_0=E_k$, and so, the 3-dimensional space $E_0 \oplus E_{\pm \omega i}$ is invariant under $A_G$ but it is not a 
polydiagonal space. That is, $E_0 \oplus E_{\pm \omega i}$ is not a synchrony space. 
\end{proof}

\ 
\vspace{3mm} 

We may ask, for each of the synchrony lattice structures in Figure~\ref{f:lista}, and for each integer $k > 1$, if there are regular  four-cell networks of valency $k$  with that lattice structure. 

For the lattice structure ${\mathcal L}_1$ in Figure~\ref{f:lista}, the following two results show that although there are no zero-Hopf valency two four-cell networks with that synchrony lattice structure, there are regular  four-cell networks of valency $2k$ for every $k$ an even positive integer with the lattice structure ${\mathcal L}_1$.

\begin{thm}\label{thm:redes_L1}
There are no zero-Hopf valency two four-cell networks with the lattice structure ${\mathcal L}_1$ in Figure~\ref{f:lista}.
\end{thm}

\begin{proof} A zero-Hopf valency two four-cell network $G$ has the lattice structure ${\mathcal L}_1$ in Figure~\ref{f:lista} if and only if it does not have any nontrivial synchrony subspace.  As $G$ is a zero-Hopf network, the  adjacency matrix $A=[a_{ij}]_{4\times4}$ of $G$ has the zero real eigenvalue. This means that, either $A$ has a zero column or has linearly dependent columns (or rows). If $A$ has a zero column, up to renumeration of the cells, it has the following block structure
$$
\left( 
\begin{array}{c|c}
B  & 0 \\
\hline 
 C & 0
\end{array}
\right),
$$
with $B$ a $3 \times 3$ matrix with row sum equal to 2.
Thus, $G$ has the nontrivial synchrony subspace where $x_1=x_2=x_3$ and, so, does not have the lattice structure ${\mathcal L}_1$ in Figure~\ref{f:lista}.

Now, consider the case where $A$ does not have a zero column and has linearly dependent columns. We have two hypotheses, either $A$ has a diagonal entry equal to $2$ and the other diagonal entries zero or it has two diagonal entries equal to $1$ and the others zero. In the first case, we can consider, without loss of generality that $a_{11}=2$ and $a_{ii}=0$, $i=2,3,4$. As the eigenvalue $2$ of the network $G$ is simple, we must have another non-zero entry in the first column, which we can assume is the second one. Moreover, we must have $a_{21}=1$, otherwise if it was  $a_{21}=2$ the first two rows of $A$ would be equal and we would have the synchrony subspace where $x_1=x_2$. Thus the matrix $A$ has one of the following two structures
$$
\left( 
\begin{array}{c|ccc}
2 & 0 & 0  & 0 \\
\hline 
1 & 0 & 1  & 0 \\
*  &  \alpha_1 & 0 & \alpha_2 \\
* & \alpha_3  & \alpha_4 & 0
\end{array}
\right),
\qquad
\left( 
\begin{array}{c|ccc}
2 & 0 & 0  & 0 \\
\hline 
1 & 0 & 0  & 1 \\
*  &  \beta_1 & 0 & \beta_2 \\
* & \beta_3  & \beta_4 & 0
\end{array}
\right),
$$
where, in both cases, the submatrix of $A$ obtained by eliminating the first row and the first column must have zero determinant and no zero column.

For the first structure, the determinant of the submatrix is $\alpha_2 \alpha_3$ and cannot be zero by the following. If $\alpha_2 =0$  then the submatrix has a zero column. If $\alpha_3 =0$ then $\alpha_1=\alpha_2=1$ and we can have:  $\alpha_4=1$, in which case the second and fourth rows are equal and we have a nontrivial synchrony subspace where $x_2=x_4$; or we have $\alpha_4=2$, in which case the matrix $A$ does not have a pair of pure imaginary eigenvalues; or $\alpha_4=0$, and then the first and four rows of $A$ are equal, implying that $x_1 = x_4$ is a nontrivial synchrony subspace.

For the second structure, the determinant of the submatrix is $\beta_1 \beta_4$ and cannot be zero by the following. If $\beta_4 =0$  then the submatrix has a zero column. If $\beta_1 =0$ then $\beta_3 = \beta_4=1$ and we can have: $\beta_2=1$ where then the third row of $A$ is equal to the second row of $A$, in which case there is the nontrivial synchrony subspace where $x_2=x_3$; or $\beta_2=2$, in which case the matrix $A$ does not have a pair of pure imaginary eigenvalues; or $\beta_2=0$ and then rows 1 and 3 of $A$ are equal and $x_1=x_3$ is  a synchrony subspace.  We conclude that when $A$ has  a diagonal entry equal to $2$ and the other diagonal entries zero, the network $G$ does not have the lattice structure ${\mathcal L}_1$ in Figure~\ref{f:lista}.

Finally, we consider the case where $A$ has two diagonal entries equal to $1$ and the others zero and no zero column. Without loss of generality, we can consider $a_{11}=a_{22}=1$ and $a_{33}=a_{44}=0$. In order for $1$ not to be an eigenvalue of $A$, we have that both the first and second columns must have (at least) one more nonzero entry, which can be $1$ or $2$. Up to permutation of cells, we have then the following hypotheses:
$$
{\tiny 
D_1 = 
\left( 
\begin{array}{cccc}
1 & 1 & 0  & 0 \\
1 & 1 & 0  & 0 \\
*  &  * & 0 & * \\
*  &  * & * & 0
\end{array}
\right),\ 
D_2 = 
\left( 
\begin{array}{cccc}
1 & 0 & *  & * \\
1 & 1 & 0  & 0 \\
*  &  1 & 0 & * \\
*  &  * & * & 0
\end{array}
\right),\ 
D_3 = 
\left( 
\begin{array}{cccc}
1 & 0 & *  & * \\
1 & 1 & 0  & 0 \\
0  &  2 & 0 & 0 \\
*  &  * & * & 0
\end{array}
\right),\ 
D_4 = 
\left( 
\begin{array}{cccc}
1 & 0 & *  & * \\
0 & 1 & *  & * \\
1  &  1 & 0 & 0 \\
*  &  * & * & 0
\end{array}
\right),\ 
}$$
$$
{\tiny 
D_5 = 
\left( 
\begin{array}{cccc}
1 & 0 & *  & * \\
0 & 1 & *  & * \\
1  &  0 & 0 & 1\\
*  &  1 & * & 0
\end{array}
\right),\ 
D_6 = 
\left( 
\begin{array}{cccc}
1 & 0 & *  & * \\
0 & 1 & *  & * \\
1  &  0 & 0 & 1 \\
0  &  2 & 0 & 0
\end{array}
\right),\ 
D_7 = 
\left( 
\begin{array}{cccc}
1 & 0 & *  & * \\
0 & 1 & *  & * \\
2  &  0 & 0 & 0 \\
0  &  2 & 0 & 0
\end{array}
\right)\, .}
$$
By inspection, we see that there are no $4\times 4$ matrices of these types with eigenvalues $2,0,\pm i$ and no nontrivial synchrony spaces. 
\end{proof}

 In the following proposition, we show the existence, for each lattice structure ${\mathcal L}_i$, with $i=1,\ldots,4$, in Figure~\ref{f:lista}, of four-cell regular zero-Hopf networks with that lattice structure.

\begin{prop} \label{prop:L234_for_k_even}
For $j=1,\ldots,4$, there are  zero-Hopf four-cell regular networks with the lattice structure ${\mathcal L}_j$ in Figure~\ref{f:lista}.
\end{prop}

\begin{proof}
For $k$ an even integer, consider the four-cell  regular networks with the following four forms for the adjacency matrices:
$$
\begin{array}{lll}
(1) \left( 
\begin{array}{cccc}
0 & 0 & 0  & 2k \\
0 & k & k  & 0 \\
k & k/2 & 0  & k /2\\
0 & k/2 & k/2 & k 
\end{array}
\right)
& \quad &
(2)
\left( 
\begin{array}{cccc}
k/2 & 0 & k/2  & 0 \\
k/2 & 0 & k/2  & 0 \\
0 & 0 & 0  & k \\
0 & k/2 & 0 & k/2 
\end{array}
\right) \\
\\
(3)
\left( 
\begin{array}{cccc}
k/2 & k/2  & 0 & 0 \\
0 & 0 & k  & 0 \\
k/2 & 0 & k/2  & 0 \\
k & 0 & 0 & 0 
\end{array}
\right)
& \quad &
(4)
\left( 
\begin{array}{cccc}
k/2 & 0 & k/2  & 0 \\
k/2 & 0 & k/2  & 0 \\
0 & 0 & 0  & k \\
k/2  & 0 & 0 & k/2 
\end{array}
\right)
\end{array}\, .
$$
Adjacency matrix (1) has eigenvalues $2k$, $0$ and $\pm i k/2$ and the corresponding network has the  synchrony lattice structure ${\mathcal L}_1$. For $j=2,3,4$, the adjacency matrix (j), has eigenvalues $k$, $0$ and $\pm i k/2$ and corresponds to a network with the synchrony  lattice structure ${\mathcal L}_j$.
\end{proof}

As an immediate consequence of Proposition~\ref{prop:L234_for_k_even}, we have that:

\begin{thm}\label{thm:infinite}
There is an infinite number of four-cell regular zero-Hopf networks.
\end{thm}

We should remark that, given any four-cell zero-Hopf regular network, we can lift it to $n$-cell zero-Hopf regular networks, with $n > 4$.

The next result provides necessary and sufficient  conditions on the eigenvalue structure of four-cell regular 
zero-Hopf networks with  the lattice structures ${\mathcal L}_i$, where $i=2,3,4$.

We start by recalling the following simple observation: if  $i\omega$, for $\omega \not=0$, is an eigenvalue of the $4 \times 4$ matrix 
$A_G$, and so,  there is $v + i w \in \C^4 \setminus \{0\}$, with $v, w \in \R^4$ such that $A_G (v + i w) = \omega i (v + iw)$,  then the 2-dimensional real eigenspace 
$E_{\pm \omega i}$ of $A_G$ is $E_{\pm \omega i} = <v, w>$. 

\begin{lemma}\label{lem_eigenvectors}
Consider a  four-cell zero-Hopf  regular network $G$. Assume $E_0 = <u>$ and $E_{i \omega}= < v + i w>$ 
where $u,v,w \in \R^4$. We have that:\\
(i) The 2-dimensional space $\Delta_1 = E_k \oplus E_0$ is a synchrony space of $G$ if and only if $u$ has two types of distinct entries. \\
(ii) The 3-dimensional space $\Delta_2 = E_k \oplus E_{\pm \omega i}$  is a synchrony space of $G$ if and only if $v + iw$  has three types of distinct entries.
\end{lemma}

\begin{proof}
(i) If $u$ has two types of distinct entries, then $E_k \oplus E_0 = <(1,1,1,1), u>$ is a 2-dimensional polydiagonal, and so, it is a network synchrony space. 

If $\Delta_1 = E_k \oplus E_0$ is a synchrony space, then it is a polydiagonal. As $\Delta_1$ is 2-dimensional, the polydiagonal has to be defined by two sets 
of equalities of coordinates, that is, 
$$
\{ x \in \R^4:\, x_{j_1} = x_{j_2}, \ x_{j_3}=x_{j_4}\},  \mbox{ or } \{ x \in \R^4:\, x_{j_1} = x_{j_2} = x_{j_3} \}, 
$$
where $j_1, j_2, j_3, j_4$ are distinct and $\{j_1, j_2, j_3, j_4\} =\{ 1,2,3,4\}$. In particular,  $u \in E_0$ has to satisfy those equalities 
(and does not have the four coordinates equal as $u \not\in E_k$). \\
(ii) Recall that if  $E_{\omega i} = <v+iw>$ then  $E_{\pm \omega i} = <v,w>$. If the vector $v+iw \in \C^4$ has three types of distinct entries, 
then $v,w \in \R^4$ have also the same three types of entries. It follows so that the space $E_k \oplus E_{\pm \omega i} = <(1,1,1,1), v, w>$ 
is a 3-dimensional polydiagonal, and so, a network synchrony space.

If $\Delta_2 = E_k \oplus E_{\pm \omega i}$ is a synchrony space, then it is a 3-dimensional polydiagonal, and so it has to be defined by one set 
of equalities of coordinates, that is, 
$$
\{ x \in \R^4:\, x_{j_1} = x_{j_2}\}, 
$$
where $j_1, j_2$ are distinct and belonging to $\{ 1,2,3,4\}$. In particular, any vector in $E_{\pm \omega i}$ has to satisfy  at least the equality $x_{j_1} = x_{j_2}$. 
Now, if the 1-dimensional complex space $E_{\omega i}$ is generated by $v + iw \in \C^4$, then $v,w \in  E_{\pm \omega i}$ and so 
satisfy at least $x_{j_1} = x_{j_2}$. Thus $v_{j_1} + i w_{j_1}  = v_{j_2} + i w_{j_2}$. Finally, if  $v + iw \in \C^4$  has less  types of distinct entries, that is, 
it has only two types of entries, and so its coordinates  satisfy two equalities, 
 then $v$ and $w$ would also both have only two types of entries. In that case, $(1,1,1,1), v,w$ would be linearly dependent, a contradiction. 
Thus, $v+iw$ satisfies only the equality $v_{j_1} + i w_{j_1}  = v_{j_2} + i w_{j_2}$. 
\end{proof}

\section{Quotients of four-cell zero-Hopf networks} \label{sec_quotients} 
If $G$ is a four-cell regular zero-Hopf network with nontrivial synchrony spaces, we describe now structural features of the possible quotient networks. 

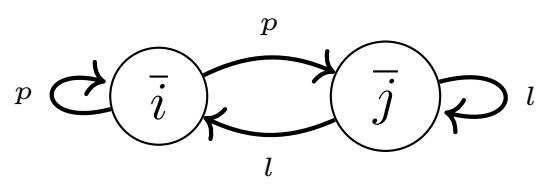
\begin{figure}
\begin{tikzpicture}
 [scale=.15,auto=left, node distance=1.5cm]
 \node[fill=white,style={circle,draw}] (n1) at (4,0) {$\overline{\small{i}}$};
  \node[fill=white,style={circle,draw}] (n2) at (14,0) {$\overline{\small{j}}$};
 \draw[->, thick] (n1) edge[thick, bend left=25]  node  [above=0.5pt] {{\tiny $p$}} (n2); 
 \draw[->, thick] (n2) edge[thick, bend left=25]  node  [below=0.5pt] {{\tiny $l$}} (n1);     
 \path
        (n1)  [->]  edge[loop left=90,thick] node {{\tiny $p$}} (n1); 
\path
        (n2)  [->]  edge[loop right=-90,thick] node {{\tiny $l$}} (n2); 
\end{tikzpicture} 
\caption{General form of a quotient network $Q_1 = G/ \Delta_1$, of a four-cell regular zero-Hopf network $G$ of valency $k \not=0$ by a two-dimensional synchrony space $\Delta_1$. Here, $p,l$ satisfy $p+l=k$ and $G$ has synchrony lattice ${\mathcal L}_3$ or ${\mathcal L}_4$ as in Figure~\ref{f:lista}. Note that $G/\Delta_1$  is ${\bf Z}_2$-symmetric if and only if $p=l=k/2$.}
\label{f:possible_2quo}
\end{figure}

\begin{prop}\label{lem_possible_2quot}
Let $G$ be a four-cell regular zero-Hopf network of valency $k > 1$ with synchrony lattice ${\mathcal L}_3$ or ${\mathcal L}_4$ as in Figure~\ref{f:lista}, and let  $\Delta_1$ be 
the  two-dimensional synchrony space. Then  the quotient network $G/\Delta_1$ 
has adjacency matrix of the form 
$$
A_{Q_1} = \left( 
\begin{array}{cc}
p & l \\
p & l
\end{array}
\right) 
$$
where $p,l$ satisfy $p + l = k$. See Figure~\ref{f:possible_2quo}. In particular, $G/\Delta_1$  is ${\bf Z}_2$-symmetric if and only if $p=l=k/2$. 
\end{prop}

\begin{proof}
By Lemma~\ref{lem_eigenvectors}\,(i), if $\Delta_1$ is a two-dimensional synchrony space of $G$, then $\Delta_1 = E_k \oplus E_0$ where $E_k = <(1,1,1,1)>$ and $E_0=<u>$ are the eigenspaces of $A_G$ associated with the valency $k$ and the zero eigenvalues, respectively, and where $u$ has two types of distinct entries. Up to re-enumeration  of cells, and scaling $u$ we may assume that $u=(1,1,1,c)$ where  $c \not= 1$ (or  $u=(0,0,0,c)$ where $c\not=0$))   or 
$u = (1,1,c,c)$ where $c \not= 1$ (or $u = (0,0,c,c)$ where $c\not=0$). 

 If $u = (1,1,1,c)$, where $c \not=1$, then $\Delta_1 = \{ x \in \R^4:\, x_1 = x_2 = x_3\}$. From $A_G u = 0$, it follows that the last column of $A_G$ has equal entries, say $l$. Take $p = k -l$. As $\Delta_1$ is a synchrony space and $A_G$ is regular, it follows that  the sum of the first three columns has all entries equal to 
$p$. The adjacency matrix of the quotient $Q_1=G/\Delta_1$ is the matrix corresponding to the restriction of $A_G$ to $\Delta_1$. Choosing the basis 
$\left( (1,1,1,0),(0,0,0,1)\right)$ of $\Delta_1$, it follows that  
$$ 
A_{Q_1} = \left( 
\begin{array}{ll}
p & l \\
p & l 
\end{array}
\right)\, .
$$
The case $u=(0,0,0,c)$ where $c\not=0$ is analogous obtaining $A_{Q_1}$ as above with $p =k$ and $l=0$. 

If $u = (1,1,c,c)$ then $\Delta_1 = \{ x \in \R^4:\, x_1 = x_2, \, x_3=x_4\}$. From $A_G u = 0$, it follows that the column vector obtained from the sum of the first two columns of $A_G$ has equal entries, say to $p$ and the column vector obtained from the sum of the last two columns of $A_G$ has also equal entries, say to $l = k -p$. Choosing the 
basis $\left( (1,1,0,0),(0,0,1,1)\right)$ of $\Delta_1$, it follows that $A_{Q_1}$ has the same form as in the previous case. The case $u = (0,0,c,c)$ where $c\not=0$ is analogous  obtaining $A_{Q_1}$ as above with $p =k$ and $l=0$. 

\end{proof}

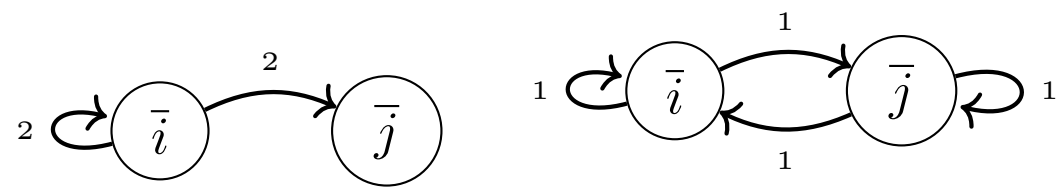
\begin{figure}
\begin{tikzpicture}
 [scale=.15,auto=left, node distance=1.5cm]
 \node[fill=white,style={circle,draw}] (n1) at (4,0) {$\overline{\small{i}}$};
  \node[fill=white,style={circle,draw}] (n2) at (14,0) {$\overline{\small{j}}$};
 \draw[->, thick] (n1) edge[thick, bend left=25]  node  [above=0.5pt] {{\tiny $2$}} (n2);  
 \path
        (n1)  [->]  edge[loop left=90,thick] node {{\tiny $2$}} (n1); 
\end{tikzpicture} \quad 
\begin{tikzpicture}
 [scale=.15,auto=left, node distance=1.5cm]
 \node[fill=white,style={circle,draw}] (n1) at (4,0) {$\overline{\small{i}}$};
  \node[fill=white,style={circle,draw}] (n2) at (14,0) {$\overline{\small{j}}$};
 \draw[->, thick] (n1) edge[thick, bend left=25]  node  [above=0.5pt] {{\tiny $1$}} (n2); 
 \draw[->, thick] (n2) edge[thick, bend left=25]  node  [below=0.5pt] {{\tiny $1$}} (n1);     
 \path
        (n1)  [->]  edge[loop left=90,thick] node {{\tiny $1$}} (n1); 
\path
        (n2)  [->]  edge[loop right=-90,thick] node {{\tiny $1$}} (n2); 
\end{tikzpicture} 
\caption{Possible two-cell quotient networks $Q_1 = G/ \Delta_1$, of a four-cell regular zero-Hopf network $G$ of valency $2$, with synchrony lattice ${\mathcal L}_3$ or ${\mathcal L}_4$ as in Figure~\ref{f:lista}, by a two-dimensional synchrony space $\Delta_1$.}
\label{f:possible_val2_int_2quo}
\end{figure}

\begin{exam} \label{example_2qval2}
If $G$ is a four-cell 
regular zero-Hopf network with valency 2
and having a two-dimensional synchrony space $\Delta_1$, then 
the possible two-cell quotient networks $Q_1 = G/\Delta_1$, up to re-enumeration of cells, have adjacency matrices 
$$
\left( 
\begin{array}{ll}
2 & 0 \\
2 & 0
\end{array}
\right) \mbox{ or }  
\left( 
\begin{array}{ll}
1 & 1 \\
1& 1 
\end{array}
\right)\, .
$$
See Figure~\ref{f:possible_val2_int_2quo}. 
\END
\end{exam}

\begin{prop}\label{prop_steady_state}
Let $G$ be a four-cell regular zero-Hopf network of valency $k > 1$ with synchrony lattice ${\mathcal L}_3$ or ${\mathcal L}_4$ as in Figure~\ref{f:lista}, and let  $\Delta_1$ be 
the  two-dimensional synchrony space with  quotient network $Q_1=G/\Delta_1$, determined by the adjacency matrix $A_{Q_1}$.
Then, generically, if $\dot{x} = F(x,\lambda)$ is a coupled cell system consistent with $G$ of the form (\ref{ccs:integers}) satisfying  $\alpha =0$, then a transcritical (resp. pitchfork) steady-state bifurcation occurs from the trivial equilibrium at $\lambda =0$ if $p\not=l$ (resp. 
if $p=l$). 
\end{prop}
\begin{proof}
It follows from Theorem 2.1. in Wang and Golubitsky~\cite{WG05} that, for the coupled systems associated with the two-cell quotient $Q_1=G/\Delta_1$ of the regular network $G$, there exists a
codimension-one synchrony-breaking bifurcation from the synchronous equilibrium  that leads
to a branch of equilibria in the corresponding synchrony subspace $\Delta_1$. The resulting bifurcation generically emerges as pitchfork or transcritical bifurcation. For two-cell networks,  pitchfork bifurcation occurs only when the network has ${\bf Z}_2$ symmetry ($p=l$ in our case). See, Golubitsky and Stewart~\cite{GS11}.
\end{proof}

\begin{prop}\label{prop_Hopf_bifurcation}
Let $G$ be a four-cell regular zero-Hopf network of valency $k>1$ with synchrony lattice ${\mathcal L}_2$ or ${\mathcal L}_4$ as in Figure~\ref{f:lista}, and let  $\Delta_2$ be 
the  three-dimensional synchrony space with  quotient network $Q_2=G/\Delta_2$. 
Then, generically, if $\dot{x} = F(x,\lambda)$ is a coupled cell system consistent with $G$ of the form (\ref{ccs:integers}) satisfying  $\alpha =0$, then a Hopf steady-state bifurcation occurs from the trivial equilibrium at $\lambda =0$. 
\end{prop}

\begin{proof}
Note that eigenvalues of $(dF)_{(0,\lambda)}$ restricted to $\Delta_2$ are 
$$\left( \frac{\partial f}{\partial u}\right)_{(0,\lambda)} + k  \left(\frac{\partial f}{\partial v_1}\right)_{(0,\lambda)} 
\quad \mbox{and} \quad \left( \frac{\partial f}{\partial u}\right)_{(0,\lambda)} \pm  i \omega   \left(\frac{\partial f}{\partial v_1}\right)_{(0,\lambda)}$$ 
in case  (\ref{ccs:integers}). Moreover, we are assuming the bifurcation condition 
 $$
 \alpha = f_u(0) = \left( \frac{\partial f}{\partial u}\right)_{(0,0)} =0\, . 
 $$
Assuming the generic hypothesis 
$$
f_{u \lambda} (0) = \left( \frac{\partial^2 f}{\partial \lambda \partial u}\right)_{(0,0)} \not=0, 
$$
the eigenvalue crossing condition of the standard Hopf Bifurcation Theorem is satisfied and so a standard Hopf bifurcation occurs from 
the origin at $\lambda =0$ leading to a unique branch of periodic solutions. 
\end{proof}

\section{Case study: four-cell valency two zero-Hopf  networks}\label{sec:val2_zHfourcells}

In this section, we enumerate, up to permutation of the cells,  the zero-Hopf  networks of four cells and valency two. We do this by considering the lattice structures  in Figure~\ref{f:lista}.
From Theorem~\ref{thm:redes_L1}, we already know that there are no four-cell zero-Hopf network with valency two with the synchrony lattice struture $\mathcal{L}_1$. So, we just need to consider the synchrony lattice strutures $\mathcal{L}_2-\mathcal{L}_4$.


In order to do the enumeration, we start by listing the possible valency two quotient networks, with two cells and three cells.

\begin{figure}
\begin{tabular}{cc}
\begin{tikzpicture}
 [scale=.15,auto=left, node distance=1.5cm, 
 ]
 \node[fill=white,style={circle,draw}] (n1) at (4,0) {$\small{1}$};
  \node[fill=white,style={circle,draw}] (n2) at (14,0) {$\small{2}$};
 \draw[->, thick] (n1) edge[thick]  node  [above=1pt, near end] {{\tiny $2$}} (n2); 
 \path
        (n1)  [->]  edge[loop left=90,thick] node {{\tiny $2$}} (n1); 
\end{tikzpicture}   \quad & \quad 
\begin{tikzpicture}
 [scale=.15,auto=left, node distance=1.5cm, 
 ]
 \node[fill=white,style={circle,draw}] (n1) at (4,0) {$\small{1}$};
  \node[fill=white,style={circle,draw}] (n2) at (14,0) {$\small{2}$};
 \draw[->, thick] (n1) edge[thick]  node  [above=1pt] {} (n2); 
 \draw[->, thick] (n2) edge[thick]  node  [above=1pt] {} (n1); 
 \path
        (n1)  [->]  edge[loop left=90,thick] node {{}} (n1); 
\path
        (n2)  [->]  edge[loop right=-90,thick] node {{}} (n2); 
\end{tikzpicture} \\
\\
$NTQ$ & $TQ$ 
\end{tabular}
\caption{The two-cell regular networks of valency two whose adjacency matrices have eigenvalues $2,0$. The network $NTQ$ is not transitive and the network $TQ$ is transitive.}
\label{f:tres}
\end{figure}
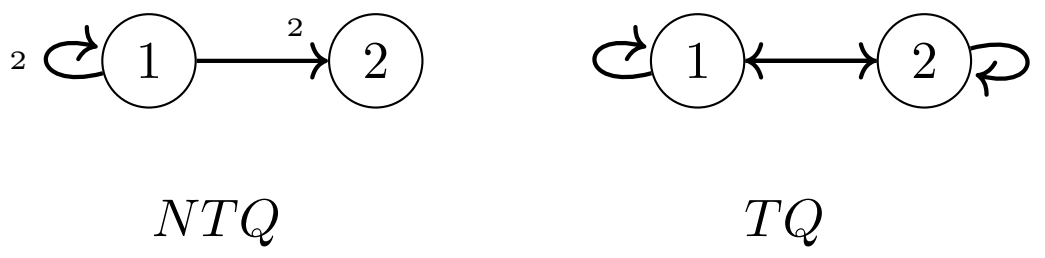

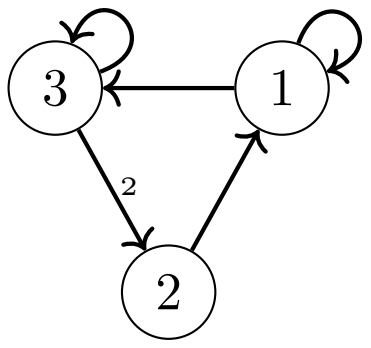
\begin{figure}
\begin{tikzpicture}
 [scale=.15,auto=left, node distance=1.5cm, 
 ]
 \node[fill=white,style={circle,draw}] (n3) at (4,0) {$\small{3}$};
  \node[fill=white,style={circle,draw}] (n1) at (14,0) {$\small{1}$};
   \node[fill=white,style={circle,draw}] (n2) at (9,-9) {$\small{2}$};
 \draw[->, thick] (n3) edge[thick]  node  [above=1pt, near end] {{\tiny $2$}} (n2); 
\draw[->, thick] (n3) to  [in=70,out=20,looseness=5] (n3);
\draw[->, thick] (n1) to  [in=20,out=70,looseness=5] (n1);
\draw[->, thick] (n1) edge[thick]  node  [above=1pt] {} (n3); 
\draw[->, thick] (n2) edge[thick]  node  [above=1pt] {} (n1); 
\end{tikzpicture}   
\caption{The three-cell regular network of valency two whose adjacency matrix has eigenvalues $2,\pm i$. It corresponds to the quotient network $Q_2$ in Figure~\ref{f:dois} and it is a transitive network.}
\label{f:quatro}
\end{figure}

\begin{prop}\label{prop:um} Up to permutation of the cells,\\
(i) there are exactly two two-cell regular  networks of valency two whose adjacency matrix has eigenvalues $2,0$:  the networks in Figure~\ref{f:tres} which have the following $2 \times 2$ adjacency matrices: 
$$
\left( 
\begin{array}{ll}
2 & 0 \\
2 & 0
\end{array}
\right),\quad 
\left( 
\begin{array}{ll}
1 & 1 \\
1 & 1
\end{array}
\right)\, .
$$
\noindent (ii) there is exactly one three-cell network  of valency two whose adjacency matrix has eigenvalues $2,\pm i$. 
See the network in Figure~\ref{f:quatro} which has the following $3 \times 3$ adjacency matrix:
$$
\left( 
\begin{array}{lll}
1 & 1 & 0 \\
0 & 0 & 2 \\
1 & 0 & 1
\end{array}
\right)\, .
$$
\end{prop}

\begin{proof}
(i) The result is a special case of Proposition~\ref{lem_possible_2quot}. Recall also Example~\ref{example_2qval2}. Equivalently, 
up to permutation of the rows and columns, there are only two $2 \times 2$ matrices with nonnegative integer entries, 
with row sum 2 and with eigenvalues $2,0$. Note that, the entries can only be $0, 1$ or $2$ and, besides row sum equal to 2, the trace of the matrix must also be equal to $2$.\\
(ii) A $3 \times 3$ matrix with the eigenvalues $2, \pm i$ has trace $2$ and determinant $2$. As the entries have to be all nonnegative integers, up to permutation of cells,  we just need to consider the following possibilities for the entries of the diagonal entries: $(2,0,0)$, or $(1,1,0)$, or $(1,0,1)$, or $(0,1,1)$. Matrices with the diagonal $(2,0,0)$ have negative or zero determinant; the other three possibilities for the diagonal entries give the following six matrices, two for each case: 
$$
\begin{array}{lll}
\left( 
\begin{array}{lll}
1 & 1 & 0 \\
0 & 1 & 1 \\
2 & 0 & 0
\end{array}
\right),\quad & \quad 
\left( 
\begin{array}{lll}
1 & 0 & 1 \\
2 & 0 & 0 \\
0 & 1 & 1
\end{array}
\right),\quad & \quad 
\left( 
\begin{array}{lll}
0 & 2 & 0 \\
0 & 1 & 1 \\
1 & 0 & 1
\end{array}
\right)\\
 \\
\left( 
\begin{array}{lll}
1 & 0 & 1 \\
1 & 1 & 0 \\
0 & 2 & 0
\end{array}
\right), \quad & \quad 
\left( 
\begin{array}{lll}
0 & 0 & 2 \\
1 & 1 & 0 \\
0 & 1 & 1
\end{array}
\right),\quad & \quad 
\left( 
\begin{array}{lll}
1 & 1 & 0 \\
0 & 0 & 2 \\
1 & 0 & 1
\end{array}
\right)\, .
\end{array}
$$
These correspond to isomorphic  networks and in Figure~\ref{f:quatro} we show a representative associated with  the last adjacency matrix. 
\end{proof}

\begin{thm}\label{thm:um}
Up to permutation of cells, there are six 
distinct zero-Hopf valency two four-cell networks that have lattice structures ${\mathcal L}_2$ or ${\mathcal L}_4$ 
in Figure~\ref{f:lista}, that is, that are lifts from a three-cell quotient with eigenvalues $2,\, \pm i$. These six classes are represented by the networks $T1, T2, T3$ in Table~\ref{t:net_syn_quo_lat1}
and the networks $ NT1, NT2, NT3$ in Table~\ref{t:net_syn_quo_lat2}. 
The networks $T1, T2, T3$  are transitive and  $NT1, NT2, NT3$ are not transitive. 
\end{thm}

\begin{proof}
From Proposition~\ref{prop:um} (ii), we have the possible valency two three-cell networks with eigenvalues $2,\, \pm i$. Using Theorem 2.5  of Aguiar, Dias, Golubistky and Leite~\cite{ADGL09}, we can describe the valency two four-cell networks whose adjacency matrix has eigenvalues $2,\, 0,\, \pm i$ which are lifts from one of these valency two three-cell networks. Moreover, the lifts will have the extra eigenvalue zero if and only if the trace of the lift is preserved. Also,  as all the networks given in Proposition~\ref{prop:um} are isomorphic, we just need to describe the lifts from one of the networks from this proposition. Choose the three-cell network $Q$ in Figure~\ref{f:quatro}  with set of cells $\{1,2,3\}$ and adjacency matrix 
$$
A_Q = 
\left( 
\begin{array}{lll}
1 & 1 & 0 \\
0 & 0 & 2 \\
1 & 0 & 1
\end{array}
\right)\, .
$$
 One of the three cells of this network must unfold into two cells. We have three hypotheses, one for each cell of $Q$. \\
(i) Unfolding the third cell of $Q$, we have to consider networks of four cells with adjacency matrix of the following block form:
$$
\left( 
\begin{array}{l|l|l}
1 & 1 & 0 \, 0 \\
\hline
0 & 0 & B_1 \\
\hline 
1 & 0 & B_2\\
1 & 0 & B_3\\
\end{array}
\right), 
$$
where $B_1$ is a block of two nonnegative integers whose sum is two and $B_2,B_3$ are both blocks of two nonnegative integers whose sum is one. We get the following six matrices which have trace two:
$$
\begin{array}{lll}
A_1 = 
\left( 
\begin{array}{cccc}
1 & 1 & 0  & 0 \\
0 & 0 & 1  & 1 \\
1 & 0 & 1 & 0\\
1 & 0 & 1 & 0
\end{array}
\right), \quad & \quad 
A_2 = 
\left( 
\begin{array}{cccc}
1 & 1 & 0  & 0 \\
0 & 0 & 1  & 1 \\
1 & 0 & 0 & 1\\
1 & 0 & 0 & 1
\end{array}
\right), \quad & \quad 
A_3 = 
\left( 
\begin{array}{cccc}
1 & 1 & 0  & 0 \\
0 & 0 & 2  & 0 \\
1 & 0 & 1 & 0\\
1 & 0 & 1 & 0
\end{array}
\right), \\
  & &  \\
A_4 = 
\left( 
\begin{array}{cccc}
1 & 1 & 0  & 0 \\
0 & 0 & 2  & 0 \\
1 & 0 & 0 & 1\\
1 & 0 & 0 & 1
\end{array}
\right), \quad & \quad 
A_5 = 
\left( 
\begin{array}{cccc}
1 & 1 & 0  & 0 \\
0 & 0 & 0  & 2 \\
1 & 0 & 1 & 0\\
1 & 0 & 1 & 0
\end{array}
\right), \quad & \quad 
A_6 = 
\left( 
\begin{array}{cccc}
1 & 1 & 0  & 0 \\
0 & 0 & 0  & 2 \\
1 & 0 & 0 & 1\\
1 & 0 & 0 & 1
\end{array}
\right)\, .
\end{array}
 $$ 
 The networks corresponding to $A_1,A_2$ are isomorphic. Also, the ones associated with $A_3,A_6$ are isomorphic, and with $A_4, A_5$ are isomorphic. \\ 
\noindent (ii) Unfolding the second cell, we get networks of four cells with adjacency matrix of the following form:
$$
\left( 
\begin{array}{l|l|l}
1 & B_1 & 0  \\
\hline
0 & 0 \, 0 & 2\\
0 & 0 \, 0 & 2\\
\hline
1 & 0 \, 0  & 1\\
\end{array}
\right), 
$$
where $B_1$ is a block of two nonnegative integers whose sum is one. We get the following matrices which have trace two:
$$
A_7 = 
\left( 
\begin{array}{cccc}
1 & 1 & 0  & 0 \\
0 & 0 & 0  & 2 \\
0 & 0 & 0 & 2\\
1 & 0 & 0 & 1
\end{array}
\right), \quad \quad 
A_8 = 
\left( 
\begin{array}{cccc}
1 & 0 & 1  & 0 \\
0 & 0 & 0  & 2 \\
0 & 0 & 0 & 2\\
1 & 0 & 0 & 1
\end{array}
\right)\, .
$$
These two matrices correspond to isomorphic networks. \\
\noindent (iii) Unfolding the first cell of $Q$, we obtain networks of four cells with adjacency matrix of the following form:
$$
\left( 
\begin{array}{l|l|l}
B_1 & 1 & 0  \\
B_2 & 1 & 0 \\
\hline
0 \, 0  & 0 & 2\\
\hline 
B_3 & 0 & 1
\end{array}
\right), 
$$
where $B_1,\, B_2,\, B_3$ are blocks of two nonnegative integers whose sum is one. We get the following matrices which have trace two: 
$$
{\tiny 
\begin{array}{llll}
A_9 = 
\left( 
\begin{array}{cccc}
1 & 0 & 1  & 0 \\
1 & 0 & 1  & 0 \\
0 & 0 & 0 & 2\\
1 & 0 & 0 & 1
\end{array}
\right), \quad & \quad 
A_{10}=
\left( 
\begin{array}{cccc}
1 & 0 & 1  & 0 \\
1 & 0 & 1  & 0 \\
0 & 0 & 0 & 2\\
0 & 1 & 0 & 1
\end{array}
\right), \ & \ 
A_{11} = 
\left( 
\begin{array}{cccc}
0 & 1 & 1  & 0 \\
0 & 1 & 1  & 0 \\
0 & 0 & 0 & 2\\
1 & 0 & 0 & 1
\end{array}
\right), \  & \ 
A_{12} = 
\left( 
\begin{array}{cccc}
0 & 1 & 1  & 0 \\
0 & 1 & 1  & 0 \\
0 & 0 & 0 & 2\\
0 & 1 & 0 & 1
\end{array}
\right)\, .
\end{array}}
$$
The networks corresponding to $A_9, \, A_{12}$ are isomorphic. Also, the networks with adjacency matrices $A_{10},\, A_{11}$ are isomorphic. 
\end{proof}

\begin{thm}\label{thm:redes_L3}
Up to permutation of cells, there are three 
distinct zero-Hopf valency two four-cell networks that have lattice structure ${\mathcal L}_3$
in Figure~\ref{f:lista}. These three classes are represented by the non-transitive networks $NT4, NT5, NT6$ in Table~\ref{t:net_syn_quo_lat3}.
\end{thm}

\begin{proof}
A network having the lattice structure ${\mathcal L}_3$ 
means that it has only one non-trivial synchrony subspace, which has dimension two, and that the eigenvalues in the restriction to that subspace are 2 and 0. 
From Proposition~\ref{prop:um} (i), up to renumeration of the cells, we have the two possible adjacency matrices of the valency two two-cell networks with eigenvalues $2,\, 0$. Using Theorem 2.5  of Aguiar, Dias, Golubistky and Leite~\cite{ADGL09}, we can describe the valency two four-cell networks whose adjacency matrix has eigenvalues $2,\, 0,\, \pm i$ which are lifts from one of these valency two two-cell networks. Moreover, by Remark 2.10 (b)-(c) in \cite{ADGL09}, if the lift is done by unfolding each of the cells into two cells, one of  the extra eigenvalues is the zero eigenvalue. Thus, the lift must be done by unfolding one of the two cells in the quotient network into three cells.  

Let us consider first the quotient network with adjacency matrix 
$$
A_Q = 
\left( 
\begin{array}{ll}
2 & 0  \\
2 & 0 
\end{array}
\right)\, .
$$
We have to unfold the first cell, otherwise the new eigenvalues of the lift adjacency matrix are zero. This way,  the last column of the lift adjacency matrix has all entries equal to $0$. Then in order to get a pair of pure imaginary eigenvalues, by Proposition~\ref{prop:um} (ii), the lift adjacency matrix must have the following structure
$$
\left( 
\begin{array}{ccc|c}
1 & 1 & 0  & 0 \\
0 & 0 & 2  & 0 \\
1 & 0 & 1 & 0\\
\hline 
 & B &  & 0
\end{array}
\right),
$$
with $B$ a $3 \times 1$ row matrix with an entry equal to $2$ and the others equal to $0$, or with two entries equal to $1$ and the other equal to $0$. From the six possible networks only the three ones in the statement of the theorem have the lattice structure ${\mathcal L}_3$, that is, have one two-dimensional synchrony subspace and no three-dimensional synchrony subspace. The three cases excluded are those where the last row of the matrix is equal to one of the other three rows which implies the existence of a three-dimensional synchrony subspace.

Now, consider the quotient network with adjacency matrix 
$$
A_Q = 
\left( 
\begin{array}{ll}
1 & 1  \\
1 & 1 
\end{array}
\right)\, .
$$
Without loss of generality, we choose to unfold the first cell. Then, the last column of the lift adjacency matrix has all entries equal to $1$ and the following structure
$$
\left( 
\begin{array}{c|c}
B & {\bf 1} \\
\hline 
C & 1
\end{array}
\right),
$$
where the rows in the $3\times 3$ matrix $B$ and in the $3 \times 1$ row matrix $C$ have one entry equal to $1$ and the others zero. It follows that, independently of the choice of $B$ and $C$, the matrix $A$  has no pair of pure imaginary eigenvalues and thus we do not get more zero-Hopf valency two four-cell networks that have the lattice structure ${\mathcal L}_3$.
\end{proof}

Tables~\ref{t:net_syn_quo_lat1}-\ref{t:net_syn_quo_lat3} list the networks representatives of the nine distinct classes of zero-Hopf valency two four-cell networks given in Theorems~\ref{thm:um}-\ref{thm:redes_L3} with the corresponding adjacency matrices eigenspaces, synchrony spaces, synchrony lattices and quotient networks.  


\section{Algebraic structure and normal forms}\label{sec_normal_form}

Normal form theory simplifies the local study of dynamical systems
through suitable near-identity coordinate transformations
\cite{SVM07}. For network dynamical systems, an additional question
is whether these transformations preserve the network architecture
and its input symmetry. We investigate this question for
\begin{equation}\label{eq:network0}
\begin{aligned}
\dot{x}_1 &= f(x_1;\overline{x_1,x_2};\lambda),\\
\dot{x}_2 &= f(x_2;\overline{x_3,x_4};\lambda),\\
\dot{x}_3 &= f(x_3;\overline{x_3,x_1};\lambda),\\
\dot{x}_4 &= f(x_4;\overline{x_1,x_3};\lambda),
\end{aligned}
\end{equation}
where the overline indicates symmetry under interchange of the
two coupling arguments:
\[
f(u;v,w;\lambda)=f(u;w,v;\lambda).
\]
The self-variable $u$ is distinguished from the coupling
variables $v,w$. Note that equations (\ref{eq:network0}) are admissible by the motivating network example in~\cite{GL09}, see Figure~\ref{f:um}. 

Our approach has three steps. We first describe the homogeneous
input-symmetric coupling functions and their associated network
vector fields. We then identify admissible generators whose Lie
brackets preserve the network structure. Finally, we use these
generators to construct a first-level normal form.

We adopt the coupled cell notation of \cite{R15}. Given maps
$\sigma_1,\ldots,\sigma_n:\{1,\ldots,N\}\to\{1,\ldots,N\}$, a
coupling function determines a network through
\begin{equation}\label{eq:coupled}
\dot{x}_j
=
f\bigl(x_{\sigma_1(j)},\ldots,x_{\sigma_n(j)}\bigr),
\qquad 1\leq j\leq N.
\end{equation}
For \eqref{eq:network0}, we choose
\begin{equation}\label{eq:sigmas}
\sigma_1=(1,2,3,4),\qquad
\sigma_2=(2,3,3,1),\qquad
\sigma_3=(1,4,1,3),
\end{equation}
where each tuple lists the values of the corresponding map.
Thus $\sigma_1$ is the identity, and the network can be written as
\begin{equation}\label{eq:networkwb}
\begin{aligned}
\dot{x}_1 &= f(x_1;x_2,x_1;\lambda),\\
\dot{x}_2 &= f(x_2;x_3,x_4;\lambda),\\
\dot{x}_3 &= f(x_3;x_3,x_1;\lambda),\\
\dot{x}_4 &= f(x_4;x_1,x_3;\lambda).
\end{aligned}
\end{equation}

Because $f$ is input-symmetric, independently reversing the two
coupling arguments in any row does not change the vector field.
There are $2^4$ choices of ordered input maps, or eight choices
up to a simultaneous interchange of $\sigma_2$ and $\sigma_3$.
These are different representations of the same input-symmetric
network class, rather than distinct vector fields.

The nontrivial input permutation is
\[
q(1)=1,\qquad q(2)=3,\qquad q(3)=2.
\]
In the terminology of \cite[Section~8]{R15}, this permutation is
a dynamical input symmetry if there exists a cell permutation
$p$ satisfying
\begin{equation}\label{eq:symmtery}
p\circ\sigma_j=\sigma_{q(j)}\circ p,
\qquad j=1,2,3.
\end{equation}
For the present network, no such permutation exists.

Indeed, any permutation satisfying \eqref{eq:symmtery} must
preserve the unordered coupling-input sets
\[
I(1)=\{1,2\},\qquad
I(2)=\{3,4\},\qquad
I(3)=I(4)=\{1,3\},
\]
in the sense that $I(p(j))=p(I(j))$.
The only two distinct cells with identical input sets are cells
$3$ and $4$, so $p$ must preserve $\{3,4\}$. Cell $3$ has a
coupling self-loop whereas cell $4$ does not, and hence both
cells are fixed. Since $I(3)=\{1,3\}$, cell $1$ is also fixed;
therefore cell $2$ is fixed as well. Thus $p$ must be the
identity. Equation \eqref{eq:symmtery} would then require
$\sigma_2=\sigma_3$, which is impossible.

Consequently, the interchange of the coupling arguments is an
input symmetry of $f$, but not a dynamical input symmetry of the
network. We therefore study admissibility-preserving
transformations directly.

\subsection{Algebraic structure}

Throughout the following analysis, the cell variables are real.
We fix the parameter at the value under consideration and
suppress $\lambda$. We assume that the origin is an equilibrium,
so that $f(0;0,0)=0$. Taylor expansions are understood formally;
equivalently, the calculations can be performed on finite Taylor
jets.

For the maps in \eqref{eq:sigmas}, define
\[
\gamma_f
=
\sum_{j=1}^{4}
f\bigl(x_{\sigma_1(j)};
       x_{\sigma_2(j)},x_{\sigma_3(j)}\bigr)\partial_j,
\qquad
\partial_j:=\frac{\partial}{\partial x_j}.
\]
We call such a vector field admissible when its coupling
function satisfies $f(u;v,w)=f(u;w,v)$.

The polynomial coupling functions belong to the invariant algebra
\[
\mathbb{R}[u,v,w]^{S_2}
=
\mathbb{R}[u,v+w,vw]
=
\mathbb{R}[u,v+w,v^2+w^2],
\]
where $S_2$ fixes $u$ and interchanges $v,w$.
The last equality follows from
\[
vw=\frac{(v+w)^2-(v^2+w^2)}{2}.
\]
Thus every input-symmetric polynomial coupling function has a
representation
\begin{equation}\label{eq:GNetwork}
f(u;v,w)=\mathfrak{f}(u,v+w,v^2+w^2).
\end{equation}
The same representation applies formally to its Taylor series.

For each $n\geq1$, introduce
\[
I_n
:=
\left\{
(i,k)\in\mathbb{Z}_{\geq0}^{\,2}:
i+k\leq n,\quad n-i-k\text{ is even}
\right\}.
\]
A basis of the homogeneous input-symmetric polynomials of
degree $n$ is
\[
\mathbf{s}^{n}_{i,k}(u;v,w)
:=
u^i(vw)^{(n-i-k)/2}(v^k+w^k),
\qquad (i,k)\in I_n.
\]
Indeed, these are the symmetrized monomials
\[
u^i\bigl(v^{r+k}w^r+v^rw^{r+k}\bigr),
\qquad r=\frac{n-i-k}{2}.
\]
For $k=0$, the expression equals $2u^i(vw)^r$; we retain this
normalization throughout.

The formal expansion of $f$ is therefore
\[
f(u;v,w)
=
\sum_{n\geq1}\sum_{(i,k)\in I_n}
a^n_{i,k}\mathbf{s}^{n}_{i,k}(u;v,w).
\]
Associate with each basis polynomial the network vector field
\begin{equation}\label{eq:basisn}
\begin{aligned}
\mathbf{S}^{n}_{i,k}
&:=\gamma_{\mathbf{s}^{n}_{i,k}}\\
&=
\sum_{j=1}^{4}
x_j^i
\bigl(x_{\sigma_2(j)}x_{\sigma_3(j)}\bigr)^{(n-i-k)/2}
\bigl(x_{\sigma_2(j)}^k+x_{\sigma_3(j)}^k\bigr)\partial_j.
\end{aligned}
\end{equation}
For network given by \eqref{eq:networkwb}, we  have
\begin{eqnarray}\label{eq:basisn-net}
\nonumber
{\bf S}^n_{i,k}&=&  {x_{{1}}}^{i} \left( x_{{1}}x_{{2}} \right)^{r} \left( {x_{
		{1}}}^{k}+{x_{{2}}}^{k} \right) \partial_1+
{x_{{2}}}^{i} \left( x_{{4}}x_{{3}} \right) {^r} \left( {x_{
		{3}}}^{k}+{x_{{4}}}^{k} \right) \partial_2
\\&&
+{x_{{3}}}^{i} \left( x_
{{3}}x_{{1}} \right)^{r} \left( {x_{{3}}}^{k}+{x_{{1}}}^{k}
\right) \partial_3+{x_{{4}}}^{i} \left( x_{{3}}x_{{1}} \right)
^{r}\left( {x_{{3}}}^{k}+{x_{{1}}}^{k} \right) \partial_4,
\end{eqnarray}
where \(\partial_i=\frac{\partial}{\partial x_i}\) and \(
r:=\frac{n-i-k}{2}.
\)
Define
\begin{equation}\label{eq:set0}
\mathscr{S}_n
:=
\operatorname{span}_{\mathbb{R}}
\left\{
\mathbf{S}^{n}_{i,k}:(i,k)\in I_n
\right\}.
\end{equation}
This is the space of homogeneous admissible vector fields of
degree $n$.

The map $f\mapsto\gamma_f$ is injective for this network:
its second component is $f(x_2;x_3,x_4)$, whose arguments are
independent variables. Consequently,
$\{\mathbf{S}^{n}_{i,k}:(i,k)\in I_n\}$ is a basis of
$\mathscr{S}_n$, and the formal network vector field has the
unique expansion
\begin{equation}\label{eq:ordg}
V^{[0]}
=
\sum_{n\geq1}\sum_{(i,k)\in I_n}
a^n_{i,k}\mathbf{S}^{n}_{i,k}.
\end{equation}

We consider the case
\[
a^1_{1,0}=0,\qquad a^1_{0,1}\neq0.
\]
Its linear part is a nonzero scalar multiple of
\[
L:=\mathbf{S}^{1}_{0,1},
\]
whose coupling function is $f_L(u;v,w)=v+w$.
Explicitly,
\[
L
=
(x_1+x_2)\partial_1
+(x_3+x_4)\partial_2
+(x_1+x_3)\partial_3
+(x_1+x_3)\partial_4.
\]
Identifying a linear vector field with its coefficient matrix,
we obtain
\[
M=
\begin{pmatrix}
1&1&0&0\\
0&0&1&1\\
1&0&1&0\\
1&0&1&0
\end{pmatrix},
\qquad
L=\sum_{j=1}^{4}(Mx)_j\partial_j.
\]
Its characteristic polynomial is
\[
\det(\mu I-M)=\mu(\mu-2)(\mu^2+1).
\]
Thus $M$ has eigenvalues $0$, $2$, and $\pm i$, with a
three-dimensional center subspace associated with $0,\pm i$.
This is the spectral configuration associated with a zero--Hopf
equilibrium; the eigenvalue calculation alone does not establish
a nondegenerate fold--Hopf bifurcation.

We carry out the normal form calculations in the full
four-dimensional phase space, without first reducing to a
center manifold; cf.\ \cite{NRS17}.

The spaces $\mathscr{S}_n$ describe the admissible homogeneous
vector fields, but this description alone does not establish
closure under Lie brackets. The next subsection examines which
generators of coordinate transformations preserve admissibility
under the bracket operations used in the normal form procedure.

\subsection{Admissible transformations}

Normal form transformations involve Lie brackets between the vector
field and the generators of coordinate changes. To retain the
network structure, we must identify generators for which these
brackets remain admissible. Although admissible vector fields form
a Lie algebra for semigroup networks \cite{R15}, this result does
not automatically imply closure for the input-symmetric class
considered here. We therefore characterize the relevant generators
directly.

Throughout this subsection, we use the Lie bracket
\[
[F,G]:=DF\,G-DG\,F.
\]
Let
\[
L:=\mathbf{S}^{1}_{0,1},
\qquad
\mathbf{M}_m
:=
\sum_{j=1}^{4}x_j^m\partial_j
=
\frac12\mathbf{S}^{m}_{m,0},
\qquad m\geq1.
\]
The factor $\frac12$ follows from the normalization
$\mathbf{s}^{m}_{m,0}=2u^m$ adopted in the preceding subsection.

Let \(g\) be an arbitrary input-symmetric homogeneous polynomial of degree \(m\). Then \(H=\gamma_g\) is an arbitrary element of \(\mathscr S_m\).
For $m\geq2$, define
\[
\mathscr{A}_m
:=
\left\{
H\in\mathscr{S}_m:
[L,H]\in\mathscr{S}_m
\right\}.
\]
Thus, $\mathscr{A}_m$ consists of the homogeneous admissible
generators whose bracket with the linear part remains admissible.
The restriction $m\geq2$ corresponds to coordinate changes tangent
to the identity.

\begin{thm}\label{thm:brackets}
The following statements hold.
\begin{enumerate}
\item
For every $m\geq2$,
\[
\mathscr{A}_m
=
\operatorname{span}_{\mathbb{R}}\{\mathbf{M}_m\}.
\]

\item
Let $f$ be an input-symmetric polynomial coupling function and
$m\geq1$. Then
\begin{equation}\label{eq:diagonal-bracket-operator}
[\gamma_f,\mathbf{M}_m]
=
\gamma_{\mathcal{D}_m f},
\qquad
\mathcal{D}_m f
:=
u^m f_u+v^m f_v+w^m f_w-mu^{m-1}f.
\end{equation}
Moreover, $\mathcal{D}_m f$ is again input-symmetric, and if $f$ is
homogeneous of degree $n$, then $\mathcal{D}_m f$ is homogeneous of
degree $n+m-1$.

\item
For every $n,m\geq1$,
\[
[\mathscr{S}_n,\mathbf{M}_m]
\subseteq
\mathscr{S}_{n+m-1}.
\]
More precisely, for every $(i,k)\in I_n$,
\begin{equation}\label{eq:admissible-bracket}
\begin{aligned}
[\mathbf{S}^{n}_{i,k},\mathbf{M}_m]
={}&
(i-m)\mathbf{S}^{n+m-1}_{i+m-1,k}\\
&+\frac{n-i+k}{2}
\mathbf{S}^{n+m-1}_{i,k+m-1}\\
&+\frac{n-i-k}{2}
\mathbf{S}^{n+m-1}_{i,|m-k-1|}.
\end{aligned}
\end{equation}
\item
In particular,
\begin{equation}\label{eq:linear-admissible-bracket}
[L,\mathbf{M}_m]
=
-m\mathbf{S}^{m}_{m-1,1}
+\mathbf{S}^{m}_{0,m}.
\end{equation}
\item
The diagonal generators satisfy
\begin{equation}\label{eq:diagonal-lie-algebra}
[\mathbf{M}_m,\mathbf{M}_{\ell}]
=
(m-\ell)\mathbf{M}_{m+\ell-1},
\qquad m,\ell\geq1.
\end{equation}
Consequently,
\[
\mathscr{M}_{\geq2}
:=
\operatorname{span}_{\mathbb{R}}
\{\mathbf{M}_m:m\geq2\}
\]
is a Lie algebra whose adjoint action preserves
$\mathscr{S}:=\bigoplus_{n\geq1}\mathscr{S}_n$.
\end{enumerate}
\end{thm}

\begin{proof}
Let $H=\gamma_g\in\mathscr{A}_m$, where $g$ is a homogeneous
input-symmetric polynomial of degree $m\geq2$. Since
$H_1=g(x_1;x_2,x_1)$ and $H_2=g(x_2;x_3,x_4)$, we have
\[
[L,H]_1
=
H_1+H_2-(x_1+x_2)\partial_1H_1
-(x_3+x_4)\partial_2H_1.
\]
Admissibility requires this expression to be independent of
$x_3,x_4$. Evaluating at $x_3=x_4=0$ therefore gives
\[
g(x_2;x_3,x_4)
=
a x_2^m+(x_3+x_4)\partial_2H_1(x_1,x_2)
\]
for some $a\in\mathbb{R}$. The left-hand side is independent
of $x_1$, so $\partial_2H_1$ is also independent of $x_1$.
By homogeneity,
\[
\partial_2H_1=c x_2^{m-1},
\qquad
g(u;v,w)=a u^m+c u^{m-1}(v+w)
\]
for some $c\in\mathbb{R}$. However, this expression for $g$
also gives
\[
\partial_2H_1=c x_1^{m-1}.
\]
Since $m\geq2$, these identities force $c=0$.
Thus $H=a\mathbf{M}_m$.

Conversely, \eqref{eq:linear-admissible-bracket} shows that
$\mathbf{M}_m\in\mathscr{A}_m$. Hence
\[
\mathscr{A}_m
=
\operatorname{span}_{\mathbb{R}}\{\mathbf{M}_m\}.
\]


We next prove the second item. Let $f$ be an input-symmetric
polynomial coupling function. The $j$-th component of $\gamma_f$ is
$f(x_j;x_{\sigma_2(j)},x_{\sigma_3(j)})$. Since
$\mathbf{M}_m=\gamma_{u^m}$, its $j$-th component is $x_j^m$, and its
Jacobian is the diagonal matrix
$\operatorname{diag}(mx_1^{m-1},\dots,mx_4^{m-1})$. The $j$-th component of
$[\gamma_f,\mathbf{M}_m]$ is therefore
\[
\begin{aligned}
&f_u(\cdot)\,x_j^m
+f_v(\cdot)\,x_{\sigma_2(j)}^m
+f_w(\cdot)\,x_{\sigma_3(j)}^m
-mx_j^{m-1}f(\cdot)\\
&\qquad=
(\mathcal{D}_m f)\bigl(x_j;x_{\sigma_2(j)},x_{\sigma_3(j)}\bigr),
\end{aligned}
\]
where $(\cdot)=(x_j;x_{\sigma_2(j)},x_{\sigma_3(j)})$. This is the
$j$-th component of $\gamma_{\mathcal{D}_m f}$, which proves
\eqref{eq:diagonal-bracket-operator}. Interchanging $v$ and $w$ in
$\mathcal{D}_m f$ interchanges the terms $v^m f_v$ and $w^m f_w$ and
leaves the remaining terms unchanged, because $f_v(u;v,w)=f_w(u;w,v)$
when $f$ is input-symmetric. Hence $\mathcal{D}_m f$ is input-symmetric.
Finally, each of $u^m f_u$, $v^m f_v$, $w^m f_w$ and $u^{m-1}f$ raises
the degree of a homogeneous $f$ by $m-1$.

The first assertion of the third item follows directly from the second
item:
\[
[\mathscr{S}_n,\mathbf{M}_m]
\subseteq\mathscr{S}_{n+m-1}.
\]
Applying $\mathcal{D}_m$ to
$\mathbf{s}^{n}_{i,k}=u^i(vw)^r(v^k+w^k)$ gives
\begin{align*}
\mathcal{D}_m\mathbf{s}^{n}_{i,k}
={}&
(i-m)u^{i+m-1}(vw)^r(v^k+w^k)\\
&+(r+k)u^i(vw)^r
\bigl(v^{k+m-1}+w^{k+m-1}\bigr)\\
&+r u^i(vw)^r
\bigl(v^kw^{m-1}+w^kv^{m-1}\bigr).
\end{align*}
The mixed term satisfies
\[
v^kw^{m-1}+w^kv^{m-1}
=
(vw)^{\min\{k,m-1\}}
\bigl(v^{|m-k-1|}+w^{|m-k-1|}\bigr).
\]
Using
\[
k+m-1
=
|m-k-1|+2\min\{k,m-1\},
\]
we can express all three terms in the basis of degree $n+m-1$.
This yields \eqref{eq:admissible-bracket}.

Finally, applying $\mathcal{D}_m$ to $v+w$ gives
\[
\mathcal{D}_m(v+w)
=
v^m+w^m-mu^{m-1}(v+w),
\]
which proves \eqref{eq:linear-admissible-bracket}.

For the last item, take $f(u;v,w)=u^\ell$, so that
$\gamma_f=\mathbf{M}_\ell$. Then
$\mathcal{D}_m u^\ell=(\ell-m)u^{m+\ell-1}$, and
\eqref{eq:diagonal-bracket-operator} gives
$[\mathbf{M}_\ell,\mathbf{M}_m]=(\ell-m)\mathbf{M}_{m+\ell-1}$,
which is \eqref{eq:diagonal-lie-algebra}. Since $m+\ell-1\geq2$
whenever $m,\ell\geq2$, the space $\mathscr{M}_{\geq2}$ is closed
under brackets. Invariance of $\mathscr{S}$ follows from the third item.
\end{proof}

\begin{rem}
The degree-one case is exceptional. Indeed,
\[
\mathscr{S}_1
=
\operatorname{span}_{\mathbb{R}}\{\mathbf{M}_1,L\},
\qquad
[L,\mathbf{M}_1]=[L,L]=0.
\]
Thus every element of $\mathscr{S}_1$ commutes with $L$, and the
uniqueness statement in Theorem~\ref{thm:brackets} applies only
to degrees $m\geq2$.

The space $\mathscr{S}$ of admissible polynomial vector fields is a
graded vector space, but it is not closed under Lie brackets.
Nevertheless, by Theorem~\ref{thm:brackets}, repeated brackets with
the generators in $\mathscr{M}_{\geq2}$ retain both the network
structure and the input symmetry. The same conclusion holds degree by
degree for formal power series, providing admissibility-preserving
generators for the normal form procedure.
\end{rem}

\subsection{Normal form study}

We now construct a first-level normal form using the
admissibility-preserving generators of the preceding subsection.
For simplicity, assume that the linear coefficient has been
normalized so that
\[
V^{[0]}=L+\sum_{m\geq2}V_m,
\qquad V_m\in\mathscr{S}_m.
\]
All normal forms below are relative to transformations generated
by $\mathscr{M}_{\geq2}$. We work formally, without asserting
convergence of the infinite normal form transformation.

For each $m\geq2$, define the restricted homological operator
\[
\mathcal{L}_m:
\operatorname{span}_{\mathbb{R}}\{\mathbf{M}_m\}
\longrightarrow\mathscr{S}_m,
\qquad
\mathcal{L}_m(H):=[L,H].
\]
Equation \eqref{eq:linear-admissible-bracket} gives
\begin{equation}\label{eq:restricted-homological-image}
\operatorname{im}\mathcal{L}_m
=
\operatorname{span}_{\mathbb{R}}
\left\{
-m\mathbf{S}^{m}_{m-1,1}
+\mathbf{S}^{m}_{0,m}
\right\}.
\end{equation}
The image is one-dimensional, so there is one scalar
normal form condition at each degree.

Two possible complementary spaces are
\begin{align}
\mathscr{N}^{\mathrm{I}}_m
&:=
\operatorname{span}_{\mathbb{R}}
\left\{
\mathbf{S}^{m}_{i,k}:
(i,k)\in I_m,\ (i,k)\neq(0,m)
\right\},
\label{eq:normal-style-I}\\
\mathscr{N}^{\mathrm{II}}_m
&:=
\operatorname{span}_{\mathbb{R}}
\left\{
\mathbf{S}^{m}_{i,k}:
(i,k)\in I_m,\ (i,k)\neq(m-1,1)
\right\}.
\label{eq:normal-style-II}
\end{align}
Both satisfy
\[
\mathscr{S}_m
=
\operatorname{im}\mathcal{L}_m
\oplus\mathscr{N}^{\mathrm{I}}_m
=
\operatorname{im}\mathcal{L}_m
\oplus\mathscr{N}^{\mathrm{II}}_m.
\]
In particular, a complement omits one of the two basis
elements occurring in
\eqref{eq:restricted-homological-image}, not both.

\begin{thm}\label{thm:first-level-normal-form}
There exists a formal coordinate transformation tangent to
the identity, obtained by composing transformations generated
by scalar multiples of $\mathbf{M}_m$, that brings $V^{[0]}$
to
\[
V^{[1]}
=
L+\sum_{m\geq2}\sum_{(i,k)\in I_m}
\widehat a^m_{i,k}\mathbf{S}^{m}_{i,k}.
\]
For each $m\geq2$, one may prescribe either
\[
\widehat a^m_{0,m}=0
\qquad\text{or}\qquad
\widehat a^m_{m-1,1}=0.
\]
These choices correspond to
$\mathscr{N}^{\mathrm{I}}_m$ and
$\mathscr{N}^{\mathrm{II}}_m$, respectively.
\end{thm}

\begin{proof}
Proceed inductively in degree. Suppose that all degrees below
$m$ have already been brought to normal form, and write the current
degree-$m$ term as
\[
V_m
=
\alpha_m\mathbf{S}^{m}_{m-1,1}
+\beta_m\mathbf{S}^{m}_{0,m}
+R_m,
\]
where $R_m$ belongs to the span of the remaining basis elements.

Choose $H_m=t_m\mathbf{M}_m$. The pushforward under its time-one
flow is represented formally by
\[
e^{\operatorname{ad}_{H_m}}V
=
V+[H_m,V]
+\frac12[H_m,[H_m,V]]+\cdots,
\qquad
\operatorname{ad}_{H_m}V:=[H_m,V].
\]
Since $H_m$ has degree $m$, its bracket with any term of degree
$j\geq2$ has degree $m+j-1>m$. The iterated bracket
$[H_m,[H_m,L]]$ also has degree $2m-1>m$.
Thus lower-degree terms remain unchanged, and the degree-$m$
term becomes
\[
V_m+[H_m,L].
\]
Using \eqref{eq:linear-admissible-bracket}, this is
\[
(\alpha_m+mt_m)\mathbf{S}^{m}_{m-1,1}
+
(\beta_m-t_m)\mathbf{S}^{m}_{0,m}
+R_m.
\]
\begin{itemize}
    \item 
For style~I, take $t_m=\beta_m$. The coefficient of
$\mathbf{S}^{m}_{0,m}$ vanishes, and the coefficient of
$\mathbf{S}^{m}_{m-1,1}$ becomes $\alpha_m+m\beta_m$.
\item
For style~II, take $t_m=-\frac{\alpha_m}{m}$. The coefficient of
$\mathbf{S}^{m}_{m-1,1}$ vanishes, and the coefficient of
$\mathbf{S}^{m}_{0,m}$ becomes $\beta_m+\frac{\alpha_m}{m}$.
\end{itemize}

All brackets remain admissible by
Theorem~\ref{thm:brackets}. Repeating this construction in
increasing degree gives the formal normal form.
\end{proof}

\begin{rem}
The choice of complementary space determines the style of the
normal form. At a given step, the two coefficients can both be
eliminated precisely when
\[
\alpha_m+m\beta_m=0.
\]
Otherwise, only one can be eliminated. The coefficients used
at each step are those of the current transformed vector field,
because lower-degree transformations may modify higher-degree
terms.
\end{rem}

\begin{exam}\label{exam:quadratic-normal-form}
Consider the quadratic input-symmetric coupling function
\[
f(u;v,w)
=
v+w+a_3(v^2+w^2)+a_4u(v+w)+a_5u^2.
\]
Its network vector field is
\[
V^{[0]}
=
L+a_3\mathbf{S}^{2}_{0,2}
+a_4\mathbf{S}^{2}_{1,1}
+a_5\mathbf{M}_2.
\]
Explicitly,
\begin{equation}\label{eq:secondexample}
\begin{aligned}
\dot x_1
&=x_1+x_2+a_3(x_1^2+x_2^2)
+a_4x_1(x_1+x_2)+a_5x_1^2,\\
\dot x_2
&=x_3+x_4+a_3(x_3^2+x_4^2)
+a_4x_2(x_3+x_4)+a_5x_2^2,\\
\dot x_3
&=x_1+x_3+a_3(x_1^2+x_3^2)
+a_4x_3(x_1+x_3)+a_5x_3^2,\\
\dot x_4
&=x_1+x_3+a_3(x_1^2+x_3^2)
+a_4x_4(x_1+x_3)+a_5x_4^2.
\end{aligned}
\end{equation}
This example omits the additional admissible quadratic
coupling term $vw$.

Take $H_2=b_5\mathbf{M}_2$. Since
\[
[\mathbf{M}_2,L]
=
2\mathbf{S}^{2}_{1,1}-\mathbf{S}^{2}_{0,2},
\]
the transformed vector field, through degree two, is
\[
\widetilde V
=
L+(a_3-b_5)\mathbf{S}^{2}_{0,2}
+(a_4+2b_5)\mathbf{S}^{2}_{1,1}
+a_5\mathbf{M}_2+O(\|x\|^3).
\]
Equivalently,
\begin{equation}\label{eq:quadratic-transformed-coupling}
\begin{aligned}
\widetilde f(u;v,w)
={}&v+w+(a_3-b_5)(v^2+w^2)\\
&+(a_4+2b_5)u(v+w)+a_5u^2
+O\bigl(\|(u,v,w)\|^3\bigr).
\end{aligned}
\end{equation}

\medskip
\noindent\textbf{Style I.}
To eliminate $v^2+w^2$, choose $b_5=a_3$. Then
\[
f^{[1]}_{\mathrm{I}}(u;v,w)
=
v+w+(2a_3+a_4)u(v+w)+a_5u^2
+O\bigl(\|(u,v,w)\|^3\bigr).
\]
Writing $A:=2a_3+a_4$, the transformed network is
\[
\begin{aligned}
\dot x_1
&=x_1+x_2+A x_1(x_1+x_2)+a_5x_1^2
+O(\|x\|^3),\\
\dot x_2
&=x_3+x_4+A x_2(x_3+x_4)+a_5x_2^2
+O(\|x\|^3),\\
\dot x_3
&=x_1+x_3+A x_3(x_1+x_3)+a_5x_3^2
+O(\|x\|^3),\\
\dot x_4
&=x_1+x_3+A x_4(x_1+x_3)+a_5x_4^2
+O(\|x\|^3).
\end{aligned}
\]

\medskip
\noindent\textbf{Style II.}
To eliminate $u(v+w)$, choose
\(
b_5=-\frac{a_4}{2}.
\)
Then
\[
f^{[1]}_{\mathrm{II}}(u;v,w)
=
v+w+\left(a_3+\frac{a_4}{2}\right)(v^2+w^2)
+a_5u^2+O\bigl(\|(u,v,w)\|^3\bigr).
\]
Writing $B:=a_3+\frac{a_4}{2}$, the transformed network is
\[
\begin{aligned}
\dot x_1
&=x_1+x_2+B(x_1^2+x_2^2)+a_5x_1^2
+O(\|x\|^3),\\
\dot x_2
&=x_3+x_4+B(x_3^2+x_4^2)+a_5x_2^2
+O(\|x\|^3),\\
\dot x_3
&=x_1+x_3+B(x_1^2+x_3^2)+a_5x_3^2
+O(\|x\|^3),\\
\dot x_4
&=x_1+x_3+B(x_1^2+x_3^2)+a_5x_4^2
+O(\|x\|^3).
\end{aligned}
\]
Although the original vector field is quadratic, its coordinate
transformation generally produces higher-order terms.
If a term $a_6vw$ is included in the original coupling function,
its coefficient is unchanged by this quadratic normal form transformation.
\end{exam}

\begin{rem}
The diagonal normal form uses only the distinguished
self-variable and the symmetry of the two coupling arguments,
not the particular input maps. For instance, for the network
\begin{equation}\label{eq:2others}
\begin{aligned}
\dot{x}_1 &= f(x_1;\overline{x_1,x_2}), \\
\dot{x}_2 &= f(x_2;\overline{x_3,x_3}),\\
\dot{x}_3 &= f(x_3;\overline{x_1,x_3}), \\
\dot{x}_4 &= f(x_4;\overline{x_1,x_1}),
\end{aligned}
\end{equation}
bracketing with $\mathbf{M}_m$ acts on $f$ through the same
operator $\mathcal{D}_m$, so the same coupling-function
normal form applies; see
Theorem~\ref{thm:universal-diagonal-normalization}. This does
not mean that such networks are dynamically equivalent, nor that
the classification in Theorem~\ref{thm:brackets} carries over.
\end{rem}

\begin{thm}\label{thm:universal-diagonal-normalization}
Let $\sigma_2,\sigma_3:\{1,\ldots,N\}\to\{1,\ldots,N\}$
be arbitrary input maps, and consider the network
\[
\dot{x}_j
=
f\bigl(x_j;x_{\sigma_2(j)},x_{\sigma_3(j)}\bigr),
\qquad j=1,\ldots,N,
\]
where $f$ is a formal coupling function satisfying
\[
f(u;v,w)=f(u;w,v),
\qquad
f(u;v,w)=v+w+O\bigl(\|(u,v,w)\|^2\bigr).
\]
Define
\[
\gamma_f
=
\sum_{j=1}^{N}
f\bigl(x_j;x_{\sigma_2(j)},x_{\sigma_3(j)}\bigr)\partial_j,
\qquad
\mathbf{M}_m
=
\sum_{j=1}^{N}x_j^m\partial_j.
\]
Then formal diagonal coordinate transformations
tangent to the identity can bring the coupling function to
\[
\widehat f
=
v+w+
\sum_{m\geq2}\sum_{(i,k)\in I_m}
\widehat a^m_{i,k}\mathbf{s}^m_{i,k},
\]
with either of the prescribed normal form conditions
\[
\widehat a^m_{0,m}=0
\qquad\text{or}\qquad
\widehat a^m_{m-1,1}=0
\]
at each degree $m\geq2$.
Here
\[
I_m
=
\left\{
(i,k)\in\mathbb{Z}_{\geq0}^{\,2}:
i+k\leq m,\quad m-i-k\text{ is even}
\right\},
\]
and
\[
\mathbf{s}^m_{i,k}
=
u^i(vw)^{(m-i-k)/2}(v^k+w^k).
\]
Thus the same coupling-function normal form is available
for every choice of $\sigma_2,\sigma_3$.
\end{thm}

\begin{proof}
The proof of Theorem~\ref{thm:brackets}(2) uses only the form
$f(x_j;x_{\sigma_2(j)},x_{\sigma_3(j)})$ of the $j$-th component of $\gamma_f$,
not the particular input maps or the number of cells, and it remains valid
when some of $x_j, x_{\sigma_2(j)}, x_{\sigma_3(j)}$ coincide: the chain rule
then adds the corresponding partial derivatives of $f$, and since the
coinciding arguments are equal, the resulting expression is unchanged. Therefore the same computation shows that
$[\gamma_f,\mathbf{M}_m]=\gamma_{\mathcal{D}_m f}$ for arbitrary input maps,
where $\mathcal{D}_m f$, defined in \eqref{eq:diagonal-bracket-operator},
is again input-symmetric.

The normal form then follows exactly as in the proof of
Theorem~\ref{thm:first-level-normal-form}, with $\mathbf{H}_m=t_m\mathbf{M}_m$;
since that calculation is carried out on the coupling function,
it does not involve the input maps.
\end{proof}

 \begin{rem}
Theorem \ref{thm:universal-diagonal-normalization} applies, in particular, to all nine networks
classified in Section~5, provided that the linear coupling
function is normalized to $v+w$. Thus, all nine networks
admit the same prescribed normal form for the coupling
function under formal diagonal coordinate transformations
tangent to the identity. The corresponding normal-form
vector fields are obtained by substituting the input maps
of each network into the transformed coupling function.
\END 
\end{rem}



\section{Conclusions}

Theorem 2.7 of Golubistky and Lauterbach~\cite{GL09} states that the center subspace of a codimension-one synchrony-breaking bifurcation in a regular network is isomorphic to a generalized eigenspace of the adjacency matrix of that network, at least when the dimension of the network nodes is greater than one. Moreover, in \cite{GL09}, it is also shown that the action of the Jacobian matrix on the center subspace is conjugate to the action of the adjacency matrix on the associated generalized eigenspace. These results suggest that the number of distinct bifurcation scenarios in regular networks may be indexed by the number of Jordan structures of adjacency matrices, rather than by the number of networks.
Even though, as pointed out in \cite{GL09},   results such as those obtained by Elmhirst and Golubistky~\cite{EG06} on nilpotent 
Hopf bifurcations have shown that more than one scenario can occur for a fixed Jordan structure, Golubistky
and Lauterbach~\cite{GL09} conjecture that each type of center subspace leads to just a small finite number of different bifurcation scenarios. 
In contrast, note that
the number of regular networks grows exponentially with the number of cells.

It was also observed in \cite{GL09} that, in general,  Theorem 2.7 does not hold if the dimension of the network nodes is one. In particular, it  was remarked that the
adjacency matrix of a regular network may have a real eigenvalue $k$ and a pair of
complex conjugate eigenvalues with real parte $k$, implying that, a bifurcation
scenario, for coupled cell systems with one-dimensional cell phase space, may lead to
a codimension-one steady-state and Hopf bifurcation occurring simultaneously for the
same bifurcation condition at a zero-Hopf singularity. Our work here was motivated by this
remark. Depending on the network structure, there can be flow-invariant synchrony spaces.
We have restricted our attention to four-cell regular networks as the simplest setting to have simultaneously a steady-state and a Hopf bifurcation that are synchrony-breaking.
We have shown that, for four-cell regular networks with any valency number,  
despite there being an infinite number of these networks that are zero-Hopf, 
there are only four possible synchrony scenarios leading to a zero-Hopf singularity, whether each of the bifurcations, steady-state and Hopf, are forced, or not, to occur in a synchrony space.

We computed the normal forms of these networks using the method developed by Rink and Sanders~\cite{R15} for coupled cell networks. The input symmetry restricts the admissible normal form transformations to those that preserve this symmetry. Consequently, fewer terms can be eliminated, and the resulting normal forms are less simplified than in the absence of input symmetry.

{\tiny 
\begin{table}[!htb]
\begin{tabular}{|l|l|l|}
\hline 
\begin{tabular}{c}
$T1$ \\ \begin{tikzpicture}
 [scale=.15,auto=left, node distance=1.5cm]
 \node[fill=white,style={circle,draw}] (n4) at (4,0) {\small{4}};
  \node[fill=white,style={circle,draw}] (n3) at (24,0) {\small{3}};
 \node[fill=white,style={circle,draw}] (n1) at (14,9)  {\small{1}};
 \node[fill=white,style={circle,draw}] (n2) at (14,-9)  {\small{2}};
 \draw[->, thick] (n1) edge[thick]  node  [above=1pt] {} (n4); 
 \draw[->, thick] (n1) edge[thick]  node  [above=1pt] {} (n3); 
 \draw[->, thick] (n4) edge[thick] node  [above=1pt] {} (n2); 
\draw[->, thick] (n3) edge[thick] node  [below=1pt]  {} (n2); 
\draw[->, thick] (n2) edge[thick] node  [below=1pt]  {} (n1); 
\draw[->, thick] (n3) edge[thick] node  [below=1pt]  {} (n4); 
\draw[->, thick] (n1) to  [in=70,out=20,looseness=5] (n1);
\draw[->, thick] (n3) to  [in=70,out=20,looseness=5] (n3);
\end{tikzpicture}  
$A_{T1} = 
\left( 
\begin{array}{cccc}
1 & 1 & 0  & 0 \\
0 & 0 & 1  & 1 \\
1 & 0 & 1 & 0\\
1 & 0 & 1 & 0
\end{array}
\right)$ \\
$\begin{array}{l}
E_2 = <(1,1,1,1)>\\
E_0 = <(1,-1,-1,1)> \\
E_{\pm i} = <(-1,0,1,1),\, (1,-2,0,0)>
\end{array}$
\end{tabular}
 &
\begin{tabular}{c}
${\mathcal L}_4$\\
\begin{tikzpicture}
 [scale=.15,auto=left, node distance=1.5cm]
\node (n3) [ellipse,draw]  at (19,9) {$\Delta_3; 2,\, 0,\, \pm i$};
\node (n2) [ellipse,draw]  at (24,2) {$\Delta_2; 2,\, \pm i$};
\node (n1) [ellipse,draw]  at (10,-3) {$\Delta_1;\, 2,\, 0$}; 
\node (n0) [ellipse,draw] at (19,-9)  {$\Delta_0;\, 2$};
\draw (n0) to (n1);
\draw (n0) to (n2);
\draw (n1) to  (n3);
\draw (n2) to  (n3);
\end{tikzpicture} \\
\\
$\begin{array}{l}
\Delta_3 = \R^4\\
\Delta_2 = \{ x:\, x_3 = x_4\}\\
\Delta_1 = \{ x:\, x_1 = x_4,\, x_2 = x_3\}\\
\Delta_0 = \{ x:\, x_1 = x_2 = x_3 = x_4\}
\end{array}$ 
\end{tabular} & 
\begin{tabular}{l}
$Q_2$  \ \begin{tikzpicture}
 [scale=.15,auto=left, node distance=1.5cm, 
 ]
 \node[fill=white,style={circle,draw}] (n3) at (4,0) {$\small{\overline{3}}$};
  \node[fill=white,style={circle,draw}] (n1) at (14,0) {$\small{\overline{1}}$};
   \node[fill=white,style={circle,draw}] (n2) at (9,-9) {$\small{\overline{2}}$};
 \draw[->, thick] (n3) edge[thick]  node  [above=1pt, near end] {{\tiny $2$}} (n2); 
\draw[->, thick] (n3) to  [in=70,out=20,looseness=5] (n3);
\draw[->, thick] (n1) to  [in=20,out=70,looseness=5] (n1);
\draw[->, thick] (n1) edge[thick]  node  [above=1pt] {} (n3); 
\draw[->, thick] (n2) edge[thick]  node  [above=1pt] {} (n1); 
\end{tikzpicture}  \\
$Q_1$ \ \begin{tikzpicture}
 [scale=.15,auto=left, node distance=1.5cm, 
 ]
 \node[fill=white,style={circle,draw}] (n1) at (4,0) {$\small{\overline{1}}$};
  \node[fill=white,style={circle,draw}] (n2) at (14,0) {$\small{\overline{2}}$};
 \draw[->, thick] (n1) edge[thick]  node  [above=1pt] {} (n2); 
 \draw[->, thick] (n2) edge[thick]  node  [above=1pt] {} (n1); 
\draw[->, thick] (n1) to  [in=70,out=20,looseness=5] (n1);
\draw[->, thick] (n2) to  [in=20,out=70,looseness=5] (n2);
\end{tikzpicture} \\
$Q_0$ \ \begin{tikzpicture}
 [scale=.15,auto=left, node distance=1.5cm]
 \node[fill=white,style={circle,draw}] (n1) at (4,0) {$\small{\overline{1}}$};
\draw[->, thick] (n1) to  [in=80,out=20,looseness=5] (n1);
\draw[->, thick] (n1) to  [in=60,out=20,looseness=5] (n1);
\end{tikzpicture} 
\end{tabular}
\\
\hline 
\begin{tabular}{c}
$T2$ \\ \begin{tikzpicture}
 [scale=.15,auto=left, node distance=1.5cm]
 \node[fill=white,style={circle,draw}] (n4) at (4,0) {\small{4}};
  \node[fill=white,style={circle,draw}] (n2) at (24,0) {\small{2}};
 \node[fill=white,style={circle,draw}] (n1) at (14,9)  {\small{1}};
 \node[fill=white,style={circle,draw}] (n3) at (14,-9)  {\small{3}};
 \draw[->, thick] (n1) edge[thick]  node  [above=1pt] {} (n4); 
 \draw[->, thick] (n1) edge[thick]  node  [above=1pt] {} (n3); 
 \draw[->, thick] (n4) edge[thick] node  [above=1pt] {} (n3); 
\draw[->, thick] (n3) edge[thick] node  [below=1pt]  {{\tiny $2$}} (n2); 
\draw[->, thick] (n2) edge[thick] node  [below=1pt]  {} (n1); 
\draw[->, thick] (n1) to  [in=70,out=20,looseness=5] (n1);
\draw[->, thick] (n4) to  [in=110,out=60,looseness=5] (n4);
\end{tikzpicture} 
$A_{T2} = 
\left( 
\begin{array}{cccc}
1 & 1 & 0  & 0 \\
0 & 0 & 2  & 0 \\
1 & 0 & 0 & 1\\
1 & 0 & 0 & 1
\end{array}
\right)$\\
$\begin{array}{l}
E_2 = <(1,1,1,1)>\\
E_0 = <(-1,1,0,1)> \\
E_{\pm i} = <(-1,0,1,1),\, (1,-2,0,0)>
\end{array}$ 
\end{tabular}
&
\begin{tabular}{c}
${\mathcal L}_2$\\
\begin{tikzpicture}
 [scale=.15,auto=left, node distance=1.5cm]
\node (n3) [ellipse,draw]  at (19,9) {$\Delta_3; 2,\, 0,\, \pm i$};
\node (n2) [ellipse,draw]  at (24,2) {$\Delta_2; 2,\, \pm i$};
\node (n0) [ellipse,draw] at (19,-9)  {$\Delta_0;\, 2$};
\draw (n0) to (n2);
\draw (n2) to (n3);
\end{tikzpicture}  \\
\\
$\begin{array}{l}
\Delta_3 = \R^4\\
\Delta_2 = \{ x:\, x_3 = x_4\}\\
\Delta_0 = \{ x:\, x_1 = x_2 = x_3 = x_4\}
\end{array}$ 
\end{tabular} & 
\begin{tabular}{l}
$Q_2$  \ \begin{tikzpicture}
 [scale=.15,auto=left, node distance=1.5cm, 
 ]
 \node[fill=white,style={circle,draw}] (n3) at (4,0) {$\small{\overline{3}}$};
  \node[fill=white,style={circle,draw}] (n1) at (14,0) {$\small{\overline{1}}$};
   \node[fill=white,style={circle,draw}] (n2) at (9,-9) {$\small{\overline{2}}$};
 \draw[->, thick] (n3) edge[thick]  node  [above=1pt, near end] {{\tiny $2$}} (n2); 
\draw[->, thick] (n3) to  [in=70,out=20,looseness=5] (n3);
\draw[->, thick] (n1) to  [in=20,out=70,looseness=5] (n1);
\draw[->, thick] (n1) edge[thick]  node  [above=1pt] {} (n3); 
\draw[->, thick] (n2) edge[thick]  node  [above=1pt] {} (n1); 
\end{tikzpicture}  \\
$Q_0$ \ \begin{tikzpicture}
 [scale=.15,auto=left, node distance=1.5cm]
 \node[fill=white,style={circle,draw}] (n1) at (4,0) {$\small{\overline{1}}$};
\draw[->, thick] (n1) to  [in=80,out=20,looseness=5] (n1);
\draw[->, thick] (n1) to  [in=60,out=20,looseness=5] (n1);
\end{tikzpicture} 
\end{tabular}\\
\hline 
\begin{tabular}{c}
$T3$ \\ \begin{tikzpicture}
 [scale=.15,auto=left, node distance=1.5cm]
 \node[fill=white,style={circle,draw}] (n2) at (4,0) {\small{2}};
  \node[fill=white,style={circle,draw}] (n3) at (24,0) {\small{3}};
 \node[fill=white,style={circle,draw}] (n1) at (14,9)  {\small{1}};
 \node[fill=white,style={circle,draw}] (n4) at (14,-9)  {\small{4}};
 \draw[->, thick] (n1) edge[thick]  node  [above=1pt] {} (n2); 
 \draw[->, thick] (n3) edge[thick] node  [above=1pt] {} (n1); 
\draw[->, thick] (n3) edge[thick] node  [below=1pt]  {} (n2); 
\draw[->, thick] (n2) edge[thick] node  [below=1pt]  {} (n4); 
\draw[->, thick] (n4) edge[thick] node  [below=1pt, near end]  {{\tiny $2$}} (n3); 
\draw[->, thick] (n3) edge[thick] node  [below=1pt]  {} (n2); 
\draw[->, thick] (n1) to  [in=70,out=20,looseness=5] (n1);
\draw[->, thick] (n4) to  [in=110,out=60,looseness=5] (n4);
\end{tikzpicture} 
$A_{T3} =
\left( 
\begin{array}{cccc}
1 & 0 & 1  & 0 \\
1 & 0 & 1  & 0 \\
0 & 0 & 0 & 2\\
0 & 1 & 0 & 1
\end{array}
\right)$ \\
$\begin{array}{l}
E_2 = <(1,1,1,1)>\\
E_0 = <(-1,0,1,0)> \\
E_{\pm i} = <(-1,-1,0,1),\, (1,1,-2,0)>
\end{array}$ 
\end{tabular} &
\begin{tabular}{c}
${\mathcal L}_2$\\
\begin{tikzpicture}
 [scale=.15,auto=left, node distance=1.5cm]
\node (n3) [ellipse,draw]  at (19,9) {$\Delta_3; 2,\, 0,\, \pm i$};
\node (n2) [ellipse,draw]  at (24,2) {$\Delta_2; 2,\, \pm i$};
\node (n0) [ellipse,draw] at (19,-9)  {$\Delta_0;\, 2$};
\draw (n0) to (n2);
\draw (n2) to (n3);
\end{tikzpicture}  \\
\\
$\begin{array}{l}
\Delta_3 = \R^4\\
\Delta_2 = \{ x:\, x_1 = x_2\}\\
\Delta_0 = \{ x:\, x_1 = x_2 = x_3 = x_4\}
\end{array}$ 
\end{tabular} & 
\begin{tabular}{l}
$Q_2$  \ \begin{tikzpicture}
 [scale=.15,auto=left, node distance=1.5cm, 
 ]
 \node[fill=white,style={circle,draw}] (n3) at (4,0) {$\small{\overline{4}}$};
  \node[fill=white,style={circle,draw}] (n1) at (14,0) {$\small{\overline{1}}$};
   \node[fill=white,style={circle,draw}] (n2) at (9,-9) {$\small{\overline{3}}$};
 \draw[->, thick] (n3) edge[thick]  node  [above=1pt, near end] {{\tiny $2$}} (n2); 
\draw[->, thick] (n3) to  [in=70,out=20,looseness=5] (n3);
\draw[->, thick] (n1) to  [in=20,out=70,looseness=5] (n1);
\draw[->, thick] (n1) edge[thick]  node  [above=1pt] {} (n3); 
\draw[->, thick] (n2) edge[thick]  node  [above=1pt] {} (n1); 
\end{tikzpicture}  \\
$Q_0$ \ \begin{tikzpicture}
 [scale=.15,auto=left, node distance=1.5cm]
 \node[fill=white,style={circle,draw}] (n1) at (4,0) {$\small{\overline{1}}$};
\draw[->, thick] (n1) to  [in=80,out=20,looseness=5] (n1);
\draw[->, thick] (n1) to  [in=60,out=20,looseness=5] (n1);
\end{tikzpicture} 
\end{tabular}\\
\hline 
\end{tabular}
\caption{Networks representatives of the three 
distinct classes of transitive zero-Hopf   valency two four-cell networks, adjacency matrices eigenspaces, synchrony spaces, synchrony lattices and quotient networks. } \label{t:net_syn_quo_lat1}
\end{table}
}

\begin{table}[!htb]
{\tiny 
\begin{tabular}{|l|l|l|}
\hline 
\begin{tabular}{c}
$NT1$ \\ \begin{tikzpicture}
 [scale=.15,auto=left, node distance=1.5cm]
 \node[fill=white,style={circle,draw}] (n4) at (4,0) {\small{4}};
  \node[fill=white,style={circle,draw}] (n2) at (24,0) {\small{2}};
 \node[fill=white,style={circle,draw}] (n1) at (14,9)  {\small{1}};
 \node[fill=white,style={circle,draw}] (n3) at (14,-9)  {\small{3}};
 \draw[->, thick] (n1) edge[thick]  node  [above=1pt] {} (n4); 
 \draw[->, thick] (n2) edge[thick]  node  [above=1pt] {} (n1); 
 \draw[->, thick] (n1) edge[thick] node  [above=1pt] {} (n3); 
\draw[->, thick] (n3) edge[thick] node  [below=1pt, near end]  {{\tiny $2$}} (n2); 
\draw[->, thick] (n3) edge[thick] node  [below=1pt]  {} (n4); 
\draw[->, thick] (n1) to  [in=70,out=20,looseness=5] (n1);
\draw[->, thick] (n3) to  [in=-40,out=10,looseness=5] (n3);
\end{tikzpicture} 
$A_{NT1} = 
\left( 
\begin{array}{cccc}
1 & 1 & 0  & 0 \\
0 & 0 & 2  & 0 \\
1 & 0 & 1 & 0\\
1 & 0 & 1 & 0
\end{array}
\right)$  \\
$\begin{array}{l}
E_2 = <(1,1,1,1)>\\
E_0 = <(0,0,0,1)> \\
E_{\pm i} = <(-1,0,1,1),\, (1,-2,0,0)>
\end{array}$
\end{tabular}
&
\begin{tabular}{c}
${\mathcal L}_4$\\
\begin{tikzpicture}
 [scale=.15,auto=left, node distance=1.5cm]
\node (n3) [ellipse,draw]  at (19,9) {$\Delta_3; 2,\, 0,\, \pm i$};
\node (n2) [ellipse,draw]  at (24,2) {$\Delta_2; 2,\, \pm i$};
\node (n1) [ellipse,draw]  at (10,-3) {$\Delta_1;\, 2,\, 0$}; 
\node (n0) [ellipse,draw] at (19,-9)  {$\Delta_0;\, 2$};
\draw (n0) to (n1);
\draw (n0) to (n2);
\draw (n1) to  (n3);
\draw (n2) to  (n3);
\end{tikzpicture} \\
\\
$\begin{array}{l}
\Delta_3 = \R^4\\
\Delta_2 = \{ x:\, x_3 = x_4\}\\
\Delta_1 = \{ x:\, x_1 = x_2 = x_3\}\\
\Delta_0 = \{ x:\, x_1 = x_2 = x_3 = x_4\}
\end{array}$ 
\end{tabular} & 
\begin{tabular}{l}
$Q_2$  \ \begin{tikzpicture}
 [scale=.15,auto=left, node distance=1.5cm, 
 ]
 \node[fill=white,style={circle,draw}] (n3) at (4,0) {$\small{\overline{3}}$};
  \node[fill=white,style={circle,draw}] (n1) at (14,0) {$\small{\overline{1}}$};
   \node[fill=white,style={circle,draw}] (n2) at (9,-9) {$\small{\overline{2}}$};
 \draw[->, thick] (n3) edge[thick]  node  [above=1pt, near end] {{\tiny $2$}} (n2); 
\draw[->, thick] (n3) to  [in=70,out=20,looseness=5] (n3);
\draw[->, thick] (n1) to  [in=20,out=70,looseness=5] (n1);
\draw[->, thick] (n1) edge[thick]  node  [above=1pt] {} (n3); 
\draw[->, thick] (n2) edge[thick]  node  [above=1pt] {} (n1); 
\end{tikzpicture}  \\
$Q_1$ \ \begin{tikzpicture}
 [scale=.15,auto=left, node distance=1.5cm, 
 ]
 \node[fill=white,style={circle,draw}] (n1) at (4,0) {$\small{\overline{1}}$};
  \node[fill=white,style={circle,draw}] (n2) at (14,0) {$\small{\overline{4}}$};
 \draw[->, thick] (n1) edge[thick]  node  [above=1pt, near end] {{\tiny $2$}} (n2); 
\draw[->, thick] (n1) to  [in=70,out=20,looseness=5] (n1);
\draw[->, thick] (n1) to  [in=90,out=20,looseness=5] (n1);
\end{tikzpicture} \\
$Q_0$ \ \begin{tikzpicture}
 [scale=.15,auto=left, node distance=1.5cm]
 \node[fill=white,style={circle,draw}] (n1) at (4,0) {$\small{\overline{1}}$};
\draw[->, thick] (n1) to  [in=80,out=20,looseness=5] (n1);
\draw[->, thick] (n1) to  [in=60,out=20,looseness=5] (n1);
\end{tikzpicture} 
\end{tabular}
\\
\hline 
\begin{tabular}{c}
$NT2$ \\ \begin{tikzpicture}
 [scale=.15,auto=left, node distance=1.5cm]
 \node[fill=white,style={circle,draw}] (n3) at (4,0) {\small{3}};
  \node[fill=white,style={circle,draw}] (n2) at (24,0) {\small{2}};
 \node[fill=white,style={circle,draw}] (n4) at (14,9)  {\small{4}};
 \node[fill=white,style={circle,draw}] (n1) at (14,-9)  {\small{1}};
 \draw[->, thick] (n1) edge[thick]  node  [above=1pt] {} (n4); 
 \draw[->, thick] (n2) edge[thick]  node  [above=1pt] {} (n1); 
 \draw[->, thick] (n4) edge[thick] node  [above=1pt, near end] {{\tiny $2$}} (n2); 
\draw[->, thick] (n4) edge[thick] node  [above=1pt, near end]  {{\tiny $2$}} (n3); 
\draw[->, thick] (n1) to  [in=-40,out=10,looseness=5] (n1);
\draw[->, thick] (n4) to  [in=70,out=20,looseness=5] (n4);
\end{tikzpicture}
$A_{NT2} = 
\left( 
\begin{array}{cccc}
1 & 1 & 0  & 0 \\
0 & 0 & 0  & 2 \\
0 & 0 & 0 & 2\\
1 & 0 & 0 & 1
\end{array}
\right)$ \\
$\begin{array}{l}
E_2 = <(1,1,1,1)>\\
E_0 = <(0,0,1,0)> \\
E_{\pm i} = <(-1,0,0,1),\, (1,-2,-2,0)>
\end{array}$ 
\end{tabular} &
\begin{tabular}{c}
${\mathcal L}_4$\\
\begin{tikzpicture}
 [scale=.15,auto=left, node distance=1.5cm]
\node (n3) [ellipse,draw]  at (19,9) {$\Delta_3; 2,\, 0,\, \pm i$};
\node (n2) [ellipse,draw]  at (24,2) {$\Delta_2; 2,\, \pm i$};
\node (n1) [ellipse,draw]  at (10,-3) {$\Delta_1;\, 2,\, 0$}; 
\node (n0) [ellipse,draw] at (19,-9)  {$\Delta_0;\, 2$};
\draw (n0) to (n1);
\draw (n0) to (n2);
\draw (n1) to  (n3);
\draw (n2) to  (n3);
\end{tikzpicture} \\
\\
$\begin{array}{l}
\Delta_3 = \R^4\\
\Delta_2 = \{ x:\, x_2 = x_3\}\\
\Delta_1 = \{ x:\, x_1 = x_2 = x_4\}\\
\Delta_0 = \{ x:\, x_1 = x_2 = x_3 = x_4\}
\end{array}$
\end{tabular} & 
\begin{tabular}{l}
$Q_2$  \ \begin{tikzpicture}
 [scale=.15,auto=left, node distance=1.5cm, 
 ]
 \node[fill=white,style={circle,draw}] (n3) at (4,0) {$\small{\overline{4}}$};
  \node[fill=white,style={circle,draw}] (n1) at (14,0) {$\small{\overline{1}}$};
   \node[fill=white,style={circle,draw}] (n2) at (9,-9) {$\small{\overline{2}}$};
 \draw[->, thick] (n3) edge[thick]  node  [above=1pt, near end] {{\tiny $2$}} (n2); 
\draw[->, thick] (n3) to  [in=70,out=20,looseness=5] (n3);
\draw[->, thick] (n1) to  [in=20,out=70,looseness=5] (n1);
\draw[->, thick] (n1) edge[thick]  node  [above=1pt] {} (n3); 
\draw[->, thick] (n2) edge[thick]  node  [above=1pt] {} (n1); 
\end{tikzpicture}  \\
$Q_1$ \ \begin{tikzpicture}
 [scale=.15,auto=left, node distance=1.5cm, 
 ]
 \node[fill=white,style={circle,draw}] (n1) at (4,0) {$\small{\overline{1}}$};
  \node[fill=white,style={circle,draw}] (n2) at (14,0) {$\small{\overline{3}}$};
 \draw[->, thick] (n1) edge[thick]  node  [above=1pt, near end] {{\tiny $2$}} (n2); 
\draw[->, thick] (n1) to  [in=70,out=20,looseness=5] (n1);
\draw[->, thick] (n1) to  [in=90,out=20,looseness=5] (n1);
\end{tikzpicture} \\
$Q_0$ \ \begin{tikzpicture}
 [scale=.15,auto=left, node distance=1.5cm]
 \node[fill=white,style={circle,draw}] (n1) at (4,0) {$\small{\overline{1}}$};
\draw[->, thick] (n1) to  [in=80,out=20,looseness=5] (n1);
\draw[->, thick] (n1) to  [in=60,out=20,looseness=5] (n1);
\end{tikzpicture} 
\end{tabular}
\\
\hline 
\begin{tabular}{c}
$NT3$ \\ \begin{tikzpicture}
 [scale=.15,auto=left, node distance=1.5cm]
 \node[fill=white,style={circle,draw}] (n1) at (4,0) {\small{1}};
  \node[fill=white,style={circle,draw}] (n3) at (24,0) {\small{3}};
 \node[fill=white,style={circle,draw}] (n4) at (14,9)  {\small{4}};
 \node[fill=white,style={circle,draw}] (n2) at (14,-9)  {\small{2}};
 \draw[->, thick] (n1) edge[thick]  node  [above=1pt] {} (n4); 
 \draw[->, thick] (n3) edge[thick]  node  [above=1pt] {} (n1); 
 \draw[->, thick] (n1) edge[thick] node  [above=1pt] {} (n2); 
\draw[->, thick] (n3) edge[thick] node  [below=1pt]  {} (n2); 
\draw[->, thick] (n4) edge[thick] node  [above=1pt, near end]  {{\tiny $2$}} (n3); 
\draw[->, thick] (n1) to  [in=100,out=50,looseness=5] (n1);
\draw[->, thick] (n4) to  [in=70,out=20,looseness=5] (n4);
\end{tikzpicture} 
$A_{NT3}= 
\left( 
\begin{array}{cccc}
1 & 0 & 1  & 0 \\
1 & 0 & 1  & 0 \\
0 & 0 & 0 & 2\\
1 & 0 & 0 & 1
\end{array}
\right)$ \\
$\begin{array}{l}
E_2 = <(1,1,1,1)>\\
E_0 = <(0,1,0,0)> \\
E_{\pm i} = <(-1,-1,0,1),\, (1,1,-2,0)>
\end{array}$ 
\end{tabular} & 
\begin{tabular}{c}
${\mathcal L}_4$\\
\begin{tikzpicture}
 [scale=.15,auto=left, node distance=1.5cm]
\node (n3) [ellipse,draw]  at (19,9) {$\Delta_3; 2,\, 0,\, \pm i$};
\node (n2) [ellipse,draw]  at (24,2) {$\Delta_2; 2,\, \pm i$};
\node (n1) [ellipse,draw]  at (10,-3) {$\Delta_1;\, 2,\, 0$}; 
\node (n0) [ellipse,draw] at (19,-9)  {$\Delta_0;\, 2$};
\draw (n0) to (n1);
\draw (n0) to (n2);
\draw (n1) to  (n3);
\draw (n2) to  (n3);
\end{tikzpicture} \\
\\
$\begin{array}{l}
\Delta_3 = \R^4\\
\Delta_2 = \{ x:\, x_1 = x_2\}\\
\Delta_1 = \{ x:\, x_1 = x_3 = x_4\}\\
\Delta_0 = \{ x:\, x_1 = x_2 = x_3 = x_4\}
\end{array}$ 
\end{tabular} & 
\begin{tabular}{l}
$Q_2$  \ \begin{tikzpicture}
 [scale=.15,auto=left, node distance=1.5cm, 
 ]
 \node[fill=white,style={circle,draw}] (n3) at (4,0) {$\small{\overline{4}}$};
  \node[fill=white,style={circle,draw}] (n1) at (14,0) {$\small{\overline{1}}$};
   \node[fill=white,style={circle,draw}] (n2) at (9,-9) {$\small{\overline{3}}$};
 \draw[->, thick] (n3) edge[thick]  node  [above=1pt, near end] {{\tiny $2$}} (n2); 
\draw[->, thick] (n3) to  [in=70,out=20,looseness=5] (n3);
\draw[->, thick] (n1) to  [in=20,out=70,looseness=5] (n1);
\draw[->, thick] (n1) edge[thick]  node  [above=1pt] {} (n3); 
\draw[->, thick] (n2) edge[thick]  node  [above=1pt] {} (n1); 
\end{tikzpicture}  \\
$Q_1$ \ \begin{tikzpicture}
 [scale=.15,auto=left, node distance=1.5cm, 
 ]
 \node[fill=white,style={circle,draw}] (n1) at (4,0) {$\small{\overline{1}}$};
  \node[fill=white,style={circle,draw}] (n2) at (14,0) {$\small{\overline{2}}$};
 \draw[->, thick] (n1) edge[thick]  node  [above=1pt, near end] {{\tiny $2$}} (n2); 
\draw[->, thick] (n1) to  [in=70,out=20,looseness=5] (n1);
\draw[->, thick] (n1) to  [in=90,out=20,looseness=5] (n1);
\end{tikzpicture} \\
$Q_0$ \ \begin{tikzpicture}
 [scale=.15,auto=left, node distance=1.5cm]
 \node[fill=white,style={circle,draw}] (n1) at (4,0) {$\small{\overline{1}}$};
\draw[->, thick] (n1) to  [in=80,out=20,looseness=5] (n1);
\draw[->, thick] (n1) to  [in=60,out=20,looseness=5] (n1);
\end{tikzpicture} 
\end{tabular}
\\
\hline 
\end{tabular}}
\caption{Networks representatives of three 
distinct classes of non transitive zero-Hopf  valency two four-cell networks, adjacency matrices eigenspaces, synchrony spaces, synchrony lattices and quotient networks.} \label{t:net_syn_quo_lat2}
\end{table}

\begin{table}[!htb]
{\tiny 
\begin{tabular}{|l|l|l|}
\hline 
\begin{tabular}{c}
$NT4$ \\ 
\begin{tikzpicture}
 [scale=.15,auto=left, node distance=1.5cm]
 \node[fill=white,style={circle,draw}] (n4) at (4,0) {\small{4}};
  \node[fill=white,style={circle,draw}] (n2) at (24,0) {\small{2}};
 \node[fill=white,style={circle,draw}] (n1) at (14,9)  {\small{1}};
 \node[fill=white,style={circle,draw}] (n3) at (14,-9)  {\small{3}};
 \draw[->, thick] (n1) edge[thick]  node  [above=1pt] {{\tiny $2$}} (n4); 
 \draw[->, thick] (n2) edge[thick]  node  [above=1pt] {} (n1); 
 \draw[->, thick] (n1) edge[thick] node  [above=1pt] {} (n3); 
\draw[->, thick] (n3) edge[thick] node  [below=1pt, near end]  {{\tiny $2$}} (n2); 
\draw[->, thick] (n1) to  [in=70,out=20,looseness=5] (n1);
\draw[->, thick] (n3) to  [in=-40,out=10,looseness=5] (n3);
\end{tikzpicture} 
$A_{NT4} = 
\left( 
\begin{array}{cccc}
1 & 1 & 0  & 0 \\
0 & 0 & 2  & 0 \\
1 & 0 & 1 & 0\\
2 & 0 & 0 & 0
\end{array}
\right)$
\\
$\begin{array}{l}
E_2 = <(1,1,1,1)>\\
E_0 = <(0,0,0,1)> \\
E_{\pm i} = <(0,-2,1,4),\, (2,-2,-1,0)>
\end{array}$
\end{tabular}
&
\begin{tabular}{c}
${\mathcal L}_3$\\
\begin{tikzpicture}
 [scale=.15,auto=left, node distance=1.5cm]
\node (n3) [ellipse,draw]  at (19,9) {$\Delta_3; 2,\, 0,\, \pm i$};
\node (n1) [ellipse,draw]  at (10,-3) {$\Delta_1;\, 2,\, 0$}; 
\node (n0) [ellipse,draw] at (19,-9)  {$\Delta_0;\, 2$};
\draw (n0) to (n1);
\draw (n1) to (n3);
\end{tikzpicture} \\
\\
$\begin{array}{l}
\Delta_3 = \R^4\\
\Delta_1 = \{ x:\, x_1 = x_2 = x_3\}\\
\Delta_0 = \{ x:\, x_1 = x_2 = x_3 = x_4\}
\end{array}$ 
\end{tabular} & 
\begin{tabular}{l}
$Q_1$ \ \begin{tikzpicture}
 [scale=.15,auto=left, node distance=1.5cm, 
 ]
 \node[fill=white,style={circle,draw}] (n1) at (4,0) {$\small{\overline{1}}$};
  \node[fill=white,style={circle,draw}] (n2) at (14,0) {$\small{\overline{4}}$};
 \draw[->, thick] (n1) edge[thick]  node  [above=1pt, near end] {{\tiny $2$}} (n2); 
\draw[->, thick] (n1) to  [in=70,out=20,looseness=5] (n1);
\draw[->, thick] (n1) to  [in=90,out=20,looseness=5] (n1);
\end{tikzpicture} \\
$Q_0$ \ \begin{tikzpicture}
 [scale=.15,auto=left, node distance=1.5cm]
 \node[fill=white,style={circle,draw}] (n1) at (4,0) {$\small{\overline{1}}$};
\draw[->, thick] (n1) to  [in=80,out=20,looseness=5] (n1);
\draw[->, thick] (n1) to  [in=60,out=20,looseness=5] (n1);
\end{tikzpicture} 
\end{tabular}
\\
\hline 
\begin{tabular}{c}
$NT5$ \\ 
\begin{tikzpicture}
 [scale=.15,auto=left, node distance=1.5cm]
 \node[fill=white,style={circle,draw}] (n4) at (4,0) {\small{4}};
  \node[fill=white,style={circle,draw}] (n2) at (24,0) {\small{2}};
 \node[fill=white,style={circle,draw}] (n1) at (14,9)  {\small{1}};
 \node[fill=white,style={circle,draw}] (n3) at (14,-9)  {\small{3}};
 \draw[->, thick] (n2) edge[thick]  node  [above=1pt] {$\ $  $\ $  {\tiny $\ 2$}} (n4); 
 \draw[->, thick] (n2) edge[thick]  node  [above=1pt] {} (n1); 
 \draw[->, thick] (n1) edge[thick] node  [above=1pt] {} (n3); 
\draw[->, thick] (n3) edge[thick] node  [below=1pt, near end]  {{\tiny $2$}} (n2); 
\draw[->, thick] (n1) to  [in=70,out=20,looseness=5] (n1);
\draw[->, thick] (n3) to  [in=-40,out=10,looseness=5] (n3);
\end{tikzpicture} 
$A_{NT5} = 
\left( 
\begin{array}{cccc}
1 & 1 & 0  & 0 \\
0 & 0 & 2  & 0 \\
1 & 0 & 1 & 0\\
0 & 2 & 0 & 0
\end{array}
\right)$
\\
$\begin{array}{l}
E_2 = <(1,1,1,1)>\\
E_0 = <(0,0,0,1)> \\
E_{\pm i} = <(1,0,-1,4),\, (-1,2,0,0)>
\end{array}$ 
\end{tabular} &
\begin{tabular}{c}
${\mathcal L}_3$\\
\begin{tikzpicture}
 [scale=.15,auto=left, node distance=1.5cm]
\node (n3) [ellipse,draw]  at (19,9) {$\Delta_3; 2,\, 0,\, \pm i$};
\node (n1) [ellipse,draw]  at (10,-3) {$\Delta_1;\, 2,\, 0$}; 
\node (n0) [ellipse,draw] at (19,-9)  {$\Delta_0;\, 2$};
\draw (n0) to (n1);
\draw (n1) to (n3);
\end{tikzpicture}\\
\\
$\begin{array}{l}
\Delta_3 = \R^4\\
\Delta_1 = \{ x:\, x_1 = x_2 = x_3\}\\
\Delta_0 = \{ x:\, x_1 = x_2 = x_3 = x_4\}
\end{array}$
\end{tabular} & 
\begin{tabular}{l}
$Q_1$ \ \begin{tikzpicture}
 [scale=.15,auto=left, node distance=1.5cm, 
 ]
 \node[fill=white,style={circle,draw}] (n1) at (4,0) {$\small{\overline{1}}$};
  \node[fill=white,style={circle,draw}] (n2) at (14,0) {$\small{\overline{4}}$};
 \draw[->, thick] (n1) edge[thick]  node  [above=1pt, near end] {{\tiny $2$}} (n2); 
\draw[->, thick] (n1) to  [in=70,out=20,looseness=5] (n1);
\draw[->, thick] (n1) to  [in=90,out=20,looseness=5] (n1);
\end{tikzpicture} \\
$Q_0$ \ \begin{tikzpicture}
 [scale=.15,auto=left, node distance=1.5cm]
 \node[fill=white,style={circle,draw}] (n1) at (4,0) {$\small{\overline{1}}$};
\draw[->, thick] (n1) to  [in=80,out=20,looseness=5] (n1);
\draw[->, thick] (n1) to  [in=60,out=20,looseness=5] (n1);
\end{tikzpicture} 
\end{tabular}
\\
\hline 
\begin{tabular}{c}
$NT6$ \\ 
\begin{tikzpicture}
 [scale=.15,auto=left, node distance=1.5cm]
 \node[fill=white,style={circle,draw}] (n4) at (4,0) {\small{4}};
  \node[fill=white,style={circle,draw}] (n2) at (24,0) {\small{2}};
 \node[fill=white,style={circle,draw}] (n1) at (14,9)  {\small{1}};
 \node[fill=white,style={circle,draw}] (n3) at (14,-9)  {\small{3}};
  \draw[->, thick] (n2) edge[thick]  node  [above=1pt] {} (n4); 
 \draw[->, thick] (n3) edge[thick]  node  [above=1pt] {} (n4); 
 \draw[->, thick] (n2) edge[thick]  node  [above=1pt] {} (n1); 
 \draw[->, thick] (n1) edge[thick] node  [above=1pt] {} (n3); 
\draw[->, thick] (n3) edge[thick] node  [below=1pt, near end]  {{\tiny $2$}} (n2); 
\draw[->, thick] (n1) to  [in=70,out=20,looseness=5] (n1);
\draw[->, thick] (n3) to  [in=-40,out=10,looseness=5] (n3);
\end{tikzpicture}  
$A_{NT6} = 
\left( 
\begin{array}{cccc}
1 & 1 & 0  & 0 \\
0 & 0 & 2  & 0 \\
1 & 0 & 1 & 0\\
0 & 1 & 1 & 0
\end{array}
\right)$
 \\
$\begin{array}{l}
E_2 = <(1,1,1,1)>\\
E_0 = <(0,0,0,1)> \\
E_{\pm i} = <(1,2,-2,5),\, (-3,4,1,0)>
\end{array}$ 
\end{tabular} & 
\begin{tabular}{c}
${\mathcal L}_3$\\
\begin{tikzpicture}
 [scale=.15,auto=left, node distance=1.5cm]
\node (n3) [ellipse,draw]  at (19,9) {$\Delta_3; 2,\, 0,\, \pm i$};
\node (n1) [ellipse,draw]  at (10,-3) {$\Delta_1;\, 2,\, 0$}; 
\node (n0) [ellipse,draw] at (19,-9)  {$\Delta_0;\, 2$};
\draw (n0) to (n1);
\draw (n1) to (n3);
\end{tikzpicture}\\
\\
$\begin{array}{l}
\Delta_3 = \R^4\\
\Delta_1 = \{ x:\, x_1 = x_2 = x_3\}\\
\Delta_0 = \{ x:\, x_1 = x_2 = x_3 = x_4\}
\end{array}$ 
\end{tabular} & 
\begin{tabular}{l}
$Q_1$ \ \begin{tikzpicture}
 [scale=.15,auto=left, node distance=1.5cm, 
 ]
 \node[fill=white,style={circle,draw}] (n1) at (4,0) {$\small{\overline{1}}$};
  \node[fill=white,style={circle,draw}] (n2) at (14,0) {$\small{\overline{4}}$};
 \draw[->, thick] (n1) edge[thick]  node  [above=1pt, near end] {{\tiny $2$}} (n2); 
\draw[->, thick] (n1) to  [in=70,out=20,looseness=5] (n1);
\draw[->, thick] (n1) to  [in=90,out=20,looseness=5] (n1);
\end{tikzpicture} \\
$Q_0$ \ \begin{tikzpicture}
 [scale=.15,auto=left, node distance=1.5cm]
 \node[fill=white,style={circle,draw}] (n1) at (4,0) {$\small{\overline{1}}$};
\draw[->, thick] (n1) to  [in=80,out=20,looseness=5] (n1);
\draw[->, thick] (n1) to  [in=60,out=20,looseness=5] (n1);
\end{tikzpicture} 
\end{tabular}
\\
\hline 
\end{tabular}}
\caption{Networks representatives of three 
distinct classes of non transitive zero-Hopf  valency two four-cell networks, adjacency matrices eigenspaces, synchrony spaces, synchrony lattices and quotient networks.} \label{t:net_syn_quo_lat3}
\end{table}


\vspace{5mm}

\noindent {\bf Acknowledgments} \\
The authors thank Bob Rink and Jan Sanders for the insightful conversations on network  normal forms. MA and AD acknowledge support by the ``BeyondTheEdge: Higher-Order Networks
and Dynamics'' project (European Union, REA Grant Agreement No. 101120085) and by 
CMUP, member of LASI, which is financed by national funds through FCT -- Funda\c c\~ao para a Ci\^encia e a Tecnologia, I.P., 
under the project UID/00144 - Centro de Matem\'atica da Universidade do Porto.

\end{document}